\documentclass[12pt]{article}

\usepackage{amsfonts}
\usepackage[all]{xy}
\usepackage{epsf} 

\newtheorem{thm}{Theorem}[section]
\newtheorem{lem}{Lemma}[section]

\newtheorem{zam}{Remark}[section]
\newtheorem{pred}{Proposition}[section]

\newtheorem{cor}{Corollary}[section]

\newcommand{\dvo}{{\it Proof. ${}$}}
\newcommand{\edvo}{\rule{6pt}{6pt}}      

\newcommand{\ZZ}{{\mathbb{Z}}}
\newcommand{\QQ}{{\mathbb{Q}}}

\newcommand{\CC}{{\mathbb{C}}}
\newcommand{\PP}{{\mathbb{P}}}
\newcommand{\FF}{{\mathbb{F}}}

\newcommand{\lz}{\textquoteleft}

\title{Generic coverings of plane with A-D-E-singularities.}
\author{V.S. Kulikov  and Vic.S. Kulikov
\thanks{Partly supported by RFFI
 (No. 99-01--01133) and INTAS-OPEN-97-2072.}}
\date{        }

\begin{document}
\maketitle

\begin{abstract}
We investigate a presentation of an algebraic surface $X$ with
A-D-E-singularities as a generic covering $f: X\to \PP^2$, i.e. a finite
morphism, having at most folds and pleats apart from singular points,
isomorphic to a projection of a surface $z^2=h(x,y)$ onto the plane $x,y$
in neighbourhoods of singular points, and the branch curve $B\subset \PP^2$
of which has only nodes and ordinary cusps except singularities originated
from the singularities of $X$. It is deemed that classics proved that a
generic projection of a non-singular surface $X\subset \PP^r$ is of such
form. In this paper this result is proved for an embedding of a surface $X$
with A-D-E-singularities,
which is a composition of the given one and a Veronese embedding.
We generalize results of the paper \cite{K}, in which  Chisini's conjecture
on the unique reconstruction of $f$ by the curve $B$ is investigated. For
this fibre products of generic coverings are studied. The main inequality
bounding the degree of a covering in the case of existence of two
nonequivalent coverings with the branch curve $B$ is obtained. This inequality is used
for the proof of the Chisini conjecture for m-canonical coverings of surfaces of
general type for $m\ge 5$.
\end{abstract}

\section*{Introduction}

Let $S \subset {\PP}^r$ be a non-singular projective surface,
$f: S \to {\PP}^2 $ be its generic projection to the plane,
$B \subset {\PP}^2 $ be the branch curve, which we call the discriminant
curve. It is deemed that classics proved (see \cite{Z} , p.104) that
(i) the map $f$ is a finite covering, which has as singularities at most
double points (folds), or singular points of cuspidal type (pleats);
(ii) with this $f^*(B) =2R+C $, where the double curve $R$
is non-singular and irreducible, and the curve $C$ is reduced;
(iii) the curve $B$ is cuspidal, i.e. has at most nodes and ordinary cusps;
over a node there lie two double points, and over a cusp --
one point of cuspidal type;
(iv) the restriction of $f$ to $R$ is of degree one.
Any finite morphism $f: S \to {\PP}^2 $ is called a {\it generic}
(or {\it simple}) {\it covering}, if it possesses the same properties as a
generic projection.
Two coverings of plane $(S_1,f_1)$ and $(S_2,f_2)$ are called equivalent, if
there is a morphism $\varphi : S_1 \to S_2 $ such that
$f_1 =f_2 \circ\varphi $.

In this paper we consider a generalization of the notion of a generic
covering to the case of surfaces with A-D-E-singularities. First of all
we want to explain why we need such a generalization. A presentation of
an algebraic variety as a finite covering of the projective space is one
of the affective ways of studying projective varieties as well as their
moduli. To compare we recall what such an approach gives in the case of
curves. For a curve $C$ of genus $g$ a generic covering
$f: C \to {\PP}^1 $ is such a covering that in every fibre there is at most
one ramification point which is a double point (or a singular point of $f$ of
type $A_1$ ). Let $B \subset {\PP}^1$ be the set of branch points, and
$d=\deg B $, i.e. $d=\sharp(B) $. Then according to the Hurwitz formula
$d =2N+2g-2$, where $N=\deg f $. If $N \ge g+1$, then any curve of genus
$g$ can be presented as a simple covering of ${\PP}^1$ of degree $N$.
The set of all simple coverings (up to equivalence)
$f: C \to {\PP}^1 $ of degree $N$ with $d$ branch points is
parametrized by a Hurwitz variety $H=H^{N,d}$. Let
${\PP}^d \setminus \Delta$
($\Delta$ -- discriminant) be the projective space parametrizing the sets
of $d$ different points of ${\PP}^1$, and let ${\it M}_g$ be the moduli
space of curves of genus $g$. There are two maps: a map
$h: H\to {\PP}^d \setminus \Delta $, sending $f$ to the set of branch points
$B \subset {\PP}^1 $, and a map $\mu: H\to {\it M}_g $, sending $f$ to the
class of curves isomorphic to $C$. Hurwitz introduced and investigated
the variety $H$ in 1891. He proved that the variety $H$ is connected, and
$h$ is a finite unramified covering. In modern functorial language $H$ was
studied also by W.Fulton in 1969. The map $\mu$ is surjective (and has fibres
of dimension $N+(N-g+1)$). This gives one of the proofs ofirreducibility of
the moduli space ${\it M}_g$.

In the case of surfaces we also can consider an analog of Hurwitz variety
$H$ of all generic coverings (up to equivalence)
$f: S \to {\PP}^2 $ of degree $N$ and with discriminant curve $B$ of degree
$d$ with given number $n$ of nodes and given number $c$ of cusps.
Let ${\PP}^{\nu} ,$ $\nu =
\frac{d(d+3)}{2} $, be a projective space parametrizing curves of degree
$d$, and  $h: H \to {\PP}^{\nu} $ be a map sending a covering
$f$ to its discriminant curve $B$. In \cite{K} a Chisini
conjecture is studied. It claims that if $B$ is the discriminant curve
of a generic covering $f$ of degree $N \ge 5$, then $f$ is uniquely up to
equivalence defined by the curve $B$. In other words, it means that
the map $h$ is injective (and, besides, $N=\deg f$ is determined by $B$).
In \cite{K} it is proved that the Chisini conjecture is true for almost
all generic coverings. In particular, it is true for generic coverings
defined by a multiple canonical class.
A construction of the moduli space of surfaces of general type uses
pluricanonical maps. As is known \cite{BPV} , if $S$ is a minimal surface of
general type, then for $m \ge 5$ the linear system $\left| mK_S\right|$
blows down only (-2)-curves and gives a birational map of $S$ to a surface
$X\subset {\PP}^r$  (the canonical model) with at most A-D-E-singularities
(in other terms, rational double points, Du Val singularities, simple
singularities of Arnol'd and etc.). This requires a generalization of the
notion of a generic covering to the case of surfaces with A-D-E-singularities.

In this paper we, firstly, generalize a classical result on singularities
of generic projections of non-singular surfaces to the case of surfaces
with A-D-E-singularities. We prove that if a surface $X\subset {\PP}^r$
has at most A-D-E-singularities, then (may be after a "twist")
for a generic projection $f: X \to {\PP}^2 $ the discriminant curve $B$ also
has at most A-D-E-singularities. It follows from a slightly more general
theorem.

\begin{thm}\label{T1}
Let $X \subset {\PP}^r $ be a surface with at most isolated singularities
of the form $z^2 =h(x,y) $ (= "double planes"), $X \to {\PP}^2 $ be the
restriction to $X$ of a generic projection
${\PP}^r \backslash L \to {\PP}^2 $ from a generic linear subspace
$L$ of dimension $r-3 $. Then

(i) $f$ is a finite covering;

(ii) at non-singular points of $X$ the covering $f$ has as singularities
at most either double points (folds), or singular points of cuspidal type
(pleats); in a neighbourhood of these points $f$ is equivalent to a
projection of a surface $x=z^2$, respectively $y=z^3+xz$, to the
plane $x,y$;

(iii) in a neighbourhood of a point $s\in Sing\;X$ the covering $f$ is
analytically equivalent to a projection of a surface $z^2=h(x,y)$ to the
plane $x,y$; from (ii) and (iii) it follows that the ramification divisor
is reduced, i.e. $f^*(B)=2R+C $, where $B=f(R)$, and $R$ and $C$ are reduced
curves;

(iv) except singular points $f(Sing\:X)$ the discriminant curve $B$ is
cuspidal;

(v) the restriction of $f$ to $R$ is of degree 1.
\end{thm}

Actually, the main difficulty in the proof of this theorem lies in the
classical case, when the surface $X$ is non-singular. Unfortunately,
authors do not know a complete (and mordern) proof of this theorem, and
it seams that such a proof does not exist. Thus, the
proof, even in the case of a non-singular surface, take interest. In this
paper we prove a weakened version of Theorem \ref{T1}, in which the initial
embedding is \lz twisted' by a Veronese embedding. This is quite enough
for the purposes described above.

Thus, the curve $B$ has, firstly, \textquoteleft the same' singularities as
the surface $X$
(and as the curve $R$), which are locally defined by the equation
$h(x,y)=0$. These singularities on $B$ we call {\it s-singularities},
in particular, {\it s-nodes} and {\it s-cusps}. Besides, there are nodes
and cusps on $B$ originated from singularities of the map $f$, which we call
{\it p-nodes} and {\it p-cusps}. There are two double points of $f$ over a
p-node, at which $f$ is defined locally as a projection of surfaces $z_1=x^2$
and $z_2 =y^2$ to the plane $x,y$.

If $S$ is a surface with A-D-E-singularities, then a covering
$f: S \to {\PP}^2 $ is called {\it generic}, if it satisfies the properties
of Theorem \ref{T1}.

\vspace{0.1cm}

Secondly, we generalize the central result of \cite{K} to the case
of surfaces with A-D-E-singularities. It is proved there that if a generic
covering $f: S \to {\PP}^2 $ of a non-singular surface $S$ with discriminant
curve $B$ is of sufficiently big degree $\deg f =N$, namely under condition
$$
N > \frac{\displaystyle 4(3\bar d+g-1)}{\displaystyle 2(3\bar d+g-1) -c} \; ,
\eqno (1)
$$
where $2\bar d =\deg B$, $g$ be the geometric genus of $B$, and $c$ be
the number of cusps, then $B$ is the discriminant curve of a unique generic
covering (the Chisini conjecture holds for $B$).

We can't expect an analogous result in the case of singular surfaces,
because for a curve $B$ of even degree with at most A-D-E-singularities
there always exists a double covering, which is generic. But if two generic
coverings with given discriminant curve $B$ are coverings of sufficiently
big degree, then they are equivalent. More exactly, we prove the following
theorem. Let there are two generic coverings
$f_1: X_1 \to {\PP}^2 $ and $f_2: X_2 \to {\PP}^2 $ of surfaces with
A-D-E-singularities and with the same discriminant curve $B \subset {\PP}^2 $.
Let $f_i^*(B) =2R_i+C_i ,\; i=1,2 $. With respect to a pair of coverings
$f_1$ and $f_2$ nodes and cusps of $B$ are partitioned into four types:
ss-, sp-, ps- and pp-nodes and cusps. For example, a sp-node $b\in B$ is a
node, which is a s-node for $f_1$ and a p-node for $f_2$. The number of
sp-nodes is denoted by $n_{sp} $. Then $ n=n_{ss}+n_{sp}+n_{ps}+n_{pp} $.
The analogous terminology is used for cusps.

\begin{thm}\label{T2}
If $f_1$ and $f_2$ are nonequivalent generic coverings,
then
$$
\deg f_2 \le \frac{\displaystyle 4(3\bar d+g_1-1)}{\displaystyle 2(3\bar d+g_1-1)
-{\iota}_1} \; ,
\eqno (2)
$$
where $g_1=p_a(R_1)$ is the arithmetic genus of the curve $R_1$, and
${\iota}_1 = 2n_{sp} + 2c_{sp} + c_{pp}$.
\end{thm}

We apply the main inequality (2) to the proof of the Chisini cojecture in
the case of generic pluricanonical coverings. Let $S$ be a minimal model of
a surface of general type. According to a theorem of Bombieri \cite{BPV} ,
if $m\ge 5 $, then the m-canonical map $\varphi_m: S\to {\PP}^{p_m-1} $,
defined by the complete linear system numerically equivalent to
$\left| mK_S\right| $, is a birational
morphism, which blows down (-2)-curves on $S$. Then the canonical model
$X= \varphi_m(S)$ has at most A-D-E-singularities. A generic projection
$f: X \to {\PP}^2 $ is called a {\it generic m-canonical covering for} $S$.
We prove the following theorem.

\begin{thm}\label{T3}
Let $S_1$ and $S_2$ be minimal models of surfaces of general type with the
same $\left( K_S^2\right)$ and $\chi(S)$, and let
$f_1: X_1 \to {\PP}^2 $, $f_2: X_2 \to {\PP}^2 $ be generic m-canonical
coverings with the same discriminant curve.
Then for $m\ge 5$ the coverings $f_1$ and $f_2$  are equivalent.
\end{thm}

Consider a subvariety ${\cal H}\subset {Hilb}\times{Gr}$,
parametrizing m-canonical coverings. Here $Hilb$ is a subscheme of the
Hilbert scheme, parametrizing numerically m-canonical embeddings
$X\subset\PP^M$ of surfaces with A-D-E-singularities and fixed
$(K^2_S)$ and $\chi (S)$, $Gr$ is the Grassmann variety of projection centres
from $\PP^M$ to $\PP^2$, and $\cal H$ consists of pairs $(X\subset\PP^M ,L)$ )
such, that a restriction to $X$ of a projection with centre $L$ is a generic
covering. By theorem \ref{T1} there is a one-to-one correspondence between
the set of irreducible (respectively, connected) components of $Hilb$ and
$\cal H$. Let $h: {\cal H}\to {\PP}^{\nu}$ be a map, taking a covering to
its discriminant curve. Denote by $\cal D$ a variety of plane curves of
degree $d$ with A-D-E-singularities, among which the number of nodes
$\ge n_p$, and the number of cusps $\ge á_p$, where $d, n_p, c_p$ are defined
by invariants of $S$ (see \S 6). By theorem \ref{T3} it follows
(cf. \cite{K} , \S 5)

\noindent {\bf Corollary.} The map, induced by $h$, from the set of irreducible
(respectively, connected) components of the variety $\cal H$ to
the set of irreducible (respectively, connected) components of the variety
$\cal D$ is injective.

The proof of the main inequality (2) in \cite{K} in the case of non-singular
surfaces runs as follows. To compare two coverings $f_1$ and $f_2$,
a normalization $X$ of the fibre product $X_1\times_{{\PP}^2}X_2 $
is considered. Let $g_i: X \to X_i ,\; i=1,2 $, be the corresponding mappings
to the factors. The preimage $g_1^{-1}(R_1)=R+C$ falls into two parts, where
$R$ is the curve mapped by $g_2$ to $R_2$,
and $C$ is the curve mapped by $g_2$ to $C_2$. If
$f_1$ and $f_2$  are nonequivalent, then the surface $X$ is irreducible,
and if $X_i$ are non-singular, then $X$ is non-singular too. The main
inequality is obtained by applying the Hodge index theorem to the pair
of divisors $R$ and $C$ on $X$. We use the same idea also in the case of
surfaces with A-D-E-singularities. For this we carry out the local analysis
of the normalization of the fibre product $X$ in the case of generic
coverings of surfaces with A-D-E-singularities.

\vspace{0.1cm}

In $\S 1$ we generalize to the case of surfaces with A-D-E-singularities
the theorem on generic projections. In $\S 2$ a local analysis
of a normalization of the fibre product $X$ is carried out.
In $\S 3$ we investigate the canonical cycle of an A-D-E-singularity,
with the help of which we compute numerical invariants of a generic
covering in $\S 4$. In $\S 5$ the main inequality (2) is proved. Finally,
in $\S 6$ the Chisini conjecture for generic m-canonical coverings
of surfaces of general type is proved.

\section{Singularities of a generic projection of a surface with
A-D-E-singalarities.}

In this section we prove Theorem \ref{T1}.

{\bf 1.1.} {\it A generic projection to} ${\PP}^3$. Let $X \subset {\PP}^r$
be a surface of degree $\deg X=N $ with at most isolated hypersurface
singularities $x_1,\ldots ,x_k$, i.e. such that the dimension of the tangent
spaces $dim\;T_{X,x_i} =3$. Denote by
${\pi}_L: {\PP}^r\setminus L \to {\PP}^{e-1}$ a projection from a linear
subspace $L$ of codimension $e$. It can be obtained as a composition of
projections with centers at points. The Theorem \ref{T1} on projections of $X$ to
the plane ($e=3$) is one of a series of theorems on generic projections
for different $e$, beginning with projections from points ($e=r$) and
finishing by projections to the line ($e=2$), i.e. Lefschetz pencils.

A classical result is that, if $r>5$ $(=2\dim X +1)$, then the projection from
a generic point gives an isomorphic embedding of $X$ into ${\PP}^{r-1}$. It
follows that, if $e\ge 6$, then the projection from a generic subspace $L$
gives an isomorphic embedding of $X$ into ${\PP}^{e-1}$. In particular,
by a generic projection the surface $X$ is embedded into ${\PP}^{5}$. When
projecting to ${\PP}^4$, $e=5$, there appears isolated singularities on
${\pi}_L(X)$, which is not difficult to describe. To prove Theorem \ref{T1}
we are
going to consider a generic projection of $X$ into ${\PP}^3$, $e=4$, and to
take advantage of the following theorem.

\begin{thm}
If  $X \subset {\PP}^r $ is a surface with at most isolated hypersurface
singularities $x_i$, then the restriction of a projection
${\pi}_L: {\PP}^r \setminus L \longrightarrow {\PP}^3 $ with the centre in a
generic subspace $L \subset {\PP}^r $ of codimension 4 gives a birational
map of $X$ onto a surface $Y \subset {\PP}^3 $, which is an isomorphism
outside the double curve $D\subset X$ not passing through the points
 $x_i,$ and $Y$ has, except the points ${\pi}_L(x_i)$, at most ordinary
singularities -- the double curve $\Delta ={\pi}_L(D)$, on which there lie
a finite number of ordinary triple points  and a finite number
of pinches. In neighbourhoods of these points in appropriate
local analytic coordinates $Y$ has normal forms as follows:
$uv=0 $ for ordinary double points,
$uvw=0 $ for ordinary triple points,
$u^2-vw^2=0 $ for pinches (or "Whitney umbrellas").
\end{thm}

The contemporary proof of this theorem one can find in the textbook
\cite{G-H}. The presence of singular points $x_i$ do not add extra
troubles: we need only to see to the centre of the projection $L$ not to
intersect the tangent spaces $T_{X,x_i} ,\; dim\;T_{X,x_i} =3$. A proof of this
theorem one can find also in \cite{M} .

We want to prove that for a generic point $\xi \in {\PP}^3$ the composition
of projections ${\pi}_L$ and ${\pi}_{\xi}: {\PP}^3\setminus \xi \to {\PP}^2 $,
i.e. the projection ${\PP}^r\setminus {\pi}_L^{-1}(\xi ) \to {\PP}^2 $
with the centre  ${\pi}_L^{-1}(\xi ) $, restricted to $X$,
$f={\pi}_{\xi}\circ {{\pi}_L}_{|_X} : X \to {\PP}^2 $, gives a covering
satisfying the properties stated in Theorem \ref{T1} .

\vspace{0.1cm}

{\bf 1.2.} {\it The disposition of lines with respect to a surface} $\PP^3$.
To describe a projection ${\pi}_{\xi}$ we need to investigate the disposition
of lines $l\subset {\PP}^3$ with respect to the surface $Y$. A line $l$
is called {\it transversal} to $Y$ at a point $y$, if it is transversal to the
tangent cone to $Y$ at this point. It means that
$(l\cdot Y)_y=1$, if $y\notin Sing\; Y ;$ $(l\cdot Y)_y=2$, if
$y\in \Delta \setminus {\Delta}_t$ and $(l\cdot Y)_y=3$, if $y\in {\Delta}_t$.
We denote by ${\Delta}_t$ and ${\Delta}_p$ the set of triple points and the
set of pinches.
If $l$ is not transversal to $Y$ at a point $y$, we say that it is tangent
to $Y$ at this point. A line $l$ is called a {\it simple tangent} to $Y$
at $y$, if $y\notin Sing\; Y $ and $(l\cdot Y)_y=2$, or if $y\in \Delta
\setminus ({\Delta}_t\cup {\Delta}_p)$ and $(l\cdot Y)_y=3$, i.e.
$(l\cdot Y_i)_y=2$ for one of two branches  $Y_i$ at the point $y$. A line
$l$ is called {\it stationary} tangent, respectively {\it simple stationary}
tangent to $Y$ at $y$, if $y\notin Sing\; Y $ and $(l\cdot Y)_y\ge 3$,
respectively $=3$. A line $l$ is called {\it stationary tangent},
respectively {\it simple stationary tangent}
to $Y$, if $l$ is transversal to $Y$ at all points, except one, at which
$l$ is stationary tangent, respectively simple stationary tangent,
and, besides the other points of intersection $l\cap  Y$ are non-singular on
$Y$. Finally, $l$ is called {\it simple bitangent}, if $l$ is transversal to
$Y$ at all points, except two of them, at which the contact is simple,
the tangent planes at them are distinct, and, besides, $l\cap Sing\; Y =\emptyset $.
We want to prove that for a generic point $\xi \in {\PP}^3$ all lines
$l\ni \xi$  are at most simple bitangents and simple stationary tangents
with respect to $Y$.

   To study the disposition of lines $l\subset {\PP}^3$ with respect to
$Y $, we consider the Grassmann variety $G=G(1,3)$ and the flag variety
${\FF} =\{ (\xi ,l)\in {\PP}^3\times G \mid \; \xi \in l \} $. There are two
projections $pr_1: {\FF} \to {\PP}^3$ and $pr_2: {\FF} \to G $, which are
${\PP}^2$- and ${\PP}^1$-bundles respectively; $dim\; {\FF} =5$,
and $dim\; G=4$. In the sequal we consider  points $\xi \in {\PP}^3$
as centres of projection ${\pi}_{\xi}: {\PP}^3\setminus \xi \to {\PP}^2 $.
The fibre ${pr}_1^{-1}(\xi ) \simeq {\PP}^2$ is mapped by the projection
${pr}_2$ isomorphically onto ${\PP}_{\xi}^2 \subset G$. For $\xi \in {\PP}^3$
there is a section $s_{\xi}: {\PP}^3\setminus \xi \to {\FF}$ of the projection
${pr}_1$, $y\longmapsto (y,\overline{\xi y}) $. Then ${\pi}_{\xi}$ coincides
with the restriction of the projection ${pr}_2$ to
$s_{\xi}({\PP}^3\setminus \xi )$.

Firstly, we consider the case, when a surface $Y$ is non-singular, and then
we describe the necessary modifications and supplements in the case, when
there is a double curve $\Delta $ and isolated singularities $s_i$ on $Y$.

Consider a filtration of the variety ${\FF}$ by subvarieties
$$
Z_k=\{(\xi ,l)\in {\FF} \mid \; (l\cdot Y)_{\xi} \ge k \}.
$$
Then
$Z_1 ={pr}_1^{-1}(Y), dim\; Z_1=4$. Over a generic point $l\in G$ the map
$\varphi ={{pr}_2}_{|Z_1}: Z_1 \to G$ is an unramified covering of degree
$N$. If there are no lines on $Y$, then $\varphi $ is a finite covering,
and $Z_2$ is the ramification divisor of the covering.

Now consider restrictions of the projection ${pr}_1$. The variety $Z_2$
is isomorphic to a projectivized tangent bundle,
$Z_2\simeq {\PP}({\Theta}_Y) $, and $Z_2\to Y$ is a ${\PP}^1$-fibre bundle,
$dim\; Z_2=3$. At a generic point $y\in Y$ there are two asymptotic
directions $l_1$ and $l_2$ in $T_{Y,y}$, for which $(l_1\cdot Y)_y$ and
$(l_2\cdot Y)_y \ge 3$. Therefore, over a generic point the restriction of
${pr}_1$ onto $Z_3$, $\psi : Z_3\to Y$, is a two-sheeted covering, the branch
curve of which $P\subset Y$ is the parabolic curve consisting of points
with coinciding asymptotic directions. Some fibres of the projection ${pr}_1$
are exceptional curves of the map $\psi $. Their images on $Y$ are points
$y$, at which the restriction of the second differential of the local
equation of $Y$ onto the tangent plane $T_{Y,y}$ vanishes. Such points $y$
are called the {\it planar points} of the surface $Y$. The curve
$H = \psi (Z_4) \subset Y$ consists of points $y$, at which at least one of
the numbers $(l_i\cdot Y)_y \ge 4$
($H$ is a curve, if the surface $Y$ is not a quadric).

\vspace{0.1cm}

{\bf 1.3.} {\it Absence of non simple stationary tangents.}
Consider a product $Y\times {\FF} \subset Y\times {\PP}^3\times G$ and
projections ${pr}'_{1}$ and ${pr}'_{2}$ onto $Y\times {\PP}^3$ and
$Y\times G$. We can consider the varieties $Z_k$ as subvarieties in
$Y\times G \subset {\PP}^3\times G$. Consider a variety
$$
I_4 = \{ (y;\xi ,l)\in Y\times {\FF} \mid \; (l\cdot Y)_y \ge 4 \}
=({pr}'_{2})^{-1}(Z_4).
$$
The projection ${pr}'_{2}={id}_Y\times {pr}_2$, as well as
${pr}_2: {\FF} \to G $, is a ${\PP}^1$-bundle. Therefore, $dim\; I_4=2$ and
$dim\; {\Sigma}_4\le 2$, where ${\Sigma}_4 =p_2(I_4)$, and $p_2$ is a
projection of $Y\times{\PP}^3\times G$ to ${\PP}^3$. Then, if
$\xi\in {\PP}^3\setminus\Sigma_4$, we have that $(l\cdot Y)_y\le 3$
for any line $l\ni\xi$ at any point $y\in Y$.

\vspace{0.1cm}

{\bf 1.4.} {\it Absence of non simple bitangents.} Consider a variety
$\Sigma_{2,3}\subset {\PP}^3$, made up of non simple bitangents, and show
that $\Sigma_{2,3} \le 2$. Consider a product
$Y\times Y\times{\FF} \subset Y\times Y\times{\PP}^3\times G$ and
subvarieties $I_{i,j}$, which are closures of
$$
I^0_{i,j}= \{ (y_1,y_2;\xi ,l)\in Y\times Y\times{\FF} :
(Y\cdot l)_{y_1} \ge i , (Y\cdot l)_{y_2} \ge j , y_1\ne y_2 \} .
$$
Denote a projection of $Y\times Y\times{\FF}$ to $Y\times Y\times G$ by
${pr}_{2}'' $, and let ${pr}_{2}''(I_{i,j}) =\tilde I_{i,j}$. The projection
${pr}_{2}''$ and its restriction to $I_{i,j},$ $I_{i,j}\to \tilde I_{i,j}$
are ${\PP}^1$-bundles.

\begin{lem}\label{L1.1}
$$
\dim\; I_{2,3} \le 2 .
$$
\end{lem}

\dvo Consider subvarieties
$$
Y\times Y\times G \supset \tilde I_{1,1} \supset \tilde I_{2,1}
\supset \tilde I_{2,2} \supset \tilde I_{2,3} \: ,
$$
and let $q_1$ be a projection onto the first factor. Obviously,
$\tilde I_{1,1}$ is an irreducible variety of dimension
$dim\;\tilde I_{1,1} =4$, birationally isomorphic to $Y\times Y $.
The projection $q_1: \tilde I_{2,1} \to Y$ is a fibration, fibers of
which are curves
$q_1^{-1}(y) \simeq C_y$, where
$$
C_y = Y\cap T_{Y,y} .
$$
The curve $C_y$ has a singularity at the point $y$, which is a node for
a generic point $y$.

Furthermore, the ristriction of the projection to
$\tilde I_{2,2},$ $q_1: \tilde I_{2,2} \to Y$, is surjective, and its fibre
over a point $y\in Y$ is
$$
q_1^{-1}(y) =\{ (y,y',l) \mid \; l\subset T_{Y,y} \hphantom{a}
\mbox{and} \hphantom{a} l \hphantom{a} \mbox{is tangent to} \hphantom{a}
'_y \hphantom{a} \mbox{at} \hphantom{a} y' \} .
$$
We want to prove that $q_1(\tilde I_{2,3})$ doesn't coincide with $Y$,
i.e. the embedding $Y\subset {\PP}^3$ possesses the following property
$(L_1):$ there exists a point $y\in Y$ such that all lines $l\subset T_{Y,y}$,
passing through $y$, have at most simple contact with $C_y\setminus \{ y\}$.
We prove this below in 1.6 (Proposition \ref{P1.2}) under the
assumption that the embedding $Y\subset {\PP}^3$ is obtained by a projection
of an embedding "improved" by a Veronese embedding $v_k, \; k\ge 2$.

Thus, $\dim\;q_1(\tilde I_{2,3}) \le 1$. A generic fibre of the map
$q_1: \tilde I_{2,3} \to Y$ is of dimension zero (it being one, $Y$
is a ruled surface and we obtain a contradiction to the property $(L_1)$ ),
therefore, $\dim\;\tilde I_{2,3} \le 1$ and, consequently,
$\dim\;I_{2,3} \le 2$.
$\hphantom{aaa} \rule{6pt}{6pt}$
\vspace{0.1cm}

Set $\Sigma_{2,3} =p_3(I_{2,3})$, where $p_3$ is a projection of
$Y\times Y\times {\PP}^3\times G$ to ${\PP}^3$. It follows from Lemma \ref{L1.1}
that $dim\;\Sigma_{2,3}\le 2$. If $\xi\notin\Sigma_{2,3} $, then any line
$l\ni\xi $, touching $Y$ at two points $y_1$ and $y_2$, has a simple
contact at these points.

\vspace{0.1cm}

{\bf 1.5.} {\it Absence of 3-tangents.} Consider a product
$Y\times Y\times Y\times{\FF} \subset Y\times Y\times Y\times{\PP}^3\times G$
and subvarieties $I_{i,j,k}$,
which are closures of
$$
I^0_{i,j,k}= \{ (y_1,y_2,y_3;\xi ,l)\in Y\times Y\times Y\times{\FF} \mid \;
(Y\cdot l)_{y_1} \ge i ,
(Y\cdot l)_{y_2} \ge j , (Y\cdot l)_{y_3} \ge k  \} ,
$$
where $y_1\ne y_2\ne y_3\ne y_1$.
Denote a projection of $Y\times Y\times Y\times{\FF}$ onto
$Y\times Y\times Y\times G$ by ${pr}_{2}^{(3)}$, and let
$\tilde I_{i,j,k} ={pr}_2^{(3)}(I_{i,j,k})$.
As above, it is clear that $dim\;\tilde I_{1,1,1}=4$, and ${pr}_{1}^{(3)}$
being a ${\PP}^1$-bundle, we have $dim\; I_{1,1,1}=5$.

\begin{lem}\label{L1.2}
$$
\dim\; I_{2,2,2} \le 2 .
$$
\end{lem}

\dvo Again consider a projection of $Y\times Y\times Y\times G$ and of its
sibvarieties
$$
Y\times Y\times Y\times G \supset \tilde I_{1,1,1} \supset \tilde I_{2,1,1}
\supset \tilde I_{2,2,1} \supset \tilde I_{2,2,2} ,
$$
to the first factor. Consider $q_1: \tilde I_{2,2,2} \to Y$. For a point
$y\in Y$ we have
$
q_1^{-1}(y) =\{ (y,y_2,y_3;l) \mid \; l\subset T_{Y,y} , \:
l \hphantom{a} \mbox{is tangent to} \hphantom{a}
'_y \hphantom{a} \mbox{at points} \hphantom{a} y_2
\hphantom{a} \mbox{and} \hphantom{a} y_3\in l \} .
$
Just as in Lemma \ref{L1.1} it is sufficient to prove that $q_1(\tilde I_{2,2,2})$
doesn't coincide with $Y$. It means that there exists a point $y\in Y$,
possessing the following property $(L_2)$: none of the lines
$l\subset T_{Y,y}$, passing through $y$, is not a bitangent, i.e. can't
touch $C_y\setminus \{ y\}$ at two different points.
We prove this below in the following 1.6 (Proposition \ref{P1.2})
under the assumption that the embedding $Y\subset {\PP}^3$ is obtained by a
projection of an embedding "improved" by a Veronese embedding $v_k$.
$\hphantom{aaa} \rule{6pt}{6pt}$
\vspace{0.1cm}

Set $\Sigma_{2,2,2}=p_4(I_{2,2,2})$, where $p_4$ is a projection of
$Y\times Y\times Y\times{\PP}^3\times G$  to ${\PP}^3$. Then
$dim\;\Sigma_{2,2,2}\le 2 $ and if $\xi\notin\Sigma_{2,2,2} $,
then any line $l\ni\xi $ touches $Y$ at most at two points.

\vspace{0.1cm}

{\bf 1.6.} {\it Embeddings with a Lefschetz property.} The properties
$(L_1)$ and $(L_2)$ in the two previous subsections mean that there exists a
point $y\in Y$, for which the projection ${\pi}_y$ of the curve
$C_y\setminus \{ y\} \subset T_{Y,y}\simeq {\PP}^2$ from the point $y$
is a Lefschetz pencil. Thus, to prove Lemmas \ref{L1.1} and \ref{L1.2} it is
necessary to prove the existence of a point $y\in Y$ possessing the
following "Lefschetz property" (L) with respect to the embedding into
${\PP}^3$. We formulate it for a surface $X$ embedded into a projective
space of any dimension.

Let $X\subset {\PP}^r$ be an embedding into the projective space. We say,
that {\it a hyperplane} $L_1\subset {\PP}^r$ {\it possesses a property} (L),
if the curve $X\cap L_1$ has at most one node, i.e. $L_1$ touches $X$ at a
unique point $x$, at which the curve $X\cap L_1$ has an ordinary quadratic
singularity. In other words, the point $[L_1]\in \check{\PP}^r$,
corresponding to $L_1$, is a non-singular point of the dual variety $X^{\vee}$.
We say that a {\it pair} $(L_1,L_3)$, where $L_3\subset L_1$ is a linear
subspace of dimension $r-3 $, {\it possesses a property} (L) , if :
$L_1$ possesses the property (L) , $x\in L_3$, and a projection of the curve
$X\cap L_1 -\to {\PP}^1$ from the centre $L_3$ is a Lefschetz pencil, i.e.
any fibre of this (rational) mapping contains one singular point, and this
point is at most quadratic (is of multiplicity 2). We say that an {\it embedding}
$X\subset {\PP}^r$ {\it possesses a property} (L) , if $\exists x\in X$, for
which $L_1 =T_{X,x}$ possesses the property (L), and $L_1$ can be added to a
pair $(L_1,L_3)$ with the property (L).

It is clear that, if a pair $(L_1,L_3)$ possesses the property (L) and
$Y\subset {\PP}^3$ is obtained from $X$ by projection from a centre
$L_4\subset L_3 ,\; dim\; L_4=r-4$, then the embedding $Y\subset {\PP}^3$
possesses the property (L).

\begin{pred}\label{P1.1}
If $S\subset {\PP}^q$ is an embedding of a non-singular surface, and
$X\subset {\PP}^{r_k}$ is its composition with the Veronese embedding $v_k$
defined by polynomials of degree k, then the embedding $X\subset {\PP}^{r_k}$
possesses the property (L).
\end{pred}

\dvo  Consider the hyperplane $L_1$ corresponding to a point
$[L_1]\in X^{\vee} \setminus Sing\; X^{\vee}$. Then the curve $C=X\cap L_1$
contains a unique singular point -- a node $x\in C$.
Let $i: C \to X$ be the embedding. Consider a projection
${\pi}_{k,x}: {\PP}^{r_k}\setminus x \to {\PP}^{r_k-1}$ from the point $x$.
To prove Proposition \ref{P1.1} it is enough to show that the image
${\pi}_{k,x}(C)$ is a non-singular curve in ${\PP}^{r_k-1}$. For then,
if $L_3'$, $dim\; L_3'=r_k-3$, is a centre of projection
${\PP}^{r_k-1}\setminus L_3' \to {\PP}^1$,
which is a Lefshetz pencil for ${\pi}_{k,x}(C)$, then, obviously, the pair
$(L_1,L_3)$, where $L_3 ={\pi_{k,x}^{-1}}(L_3')\cap L_1$, possesses
the property $(L)$.

Let $I_x$ be the ideal sheaf of the point $x$ on $S$, and
$\mathcal{O}_S(1)$ be the sheaf of hyperplane sections. Under the identification
$v_k: S \simeq X$, the map ${\pi}_{k,x}$ is given by sections of
$H^0(S,\mathcal{O}_S(k)\otimes I_x)$.
Let $k=2$ and let $\sigma :\overline S \to S$ be a $\sigma$-process with
centre at the point $x$.
We can assume that $\overline S$ is embedded into $\PP^{r_2-1}$, where
$r_2-1={q(q+3)/2}$, and
the rational map $\sigma ^{-1} : S\to \overline S$ is given by sections of
$H^0(S,\mathcal{O}_S(2)\otimes I_x)$, i.e. it coincides with ${\pi}_{2,x}$.
Since the proper transform $\overline C=\sigma ^{-1}(C)\subset \overline S$
is a non-singular curve, we obtain Proposition \ref{P1.1} in the case
$k=2$. Besides, note that sections of
$i^*(H^0(S,\mathcal{O}_S(2)\otimes I_x))$ give an embedding of $\overline C$
into ${\PP}^{r_2-1}$. Consequently, for $k>2$ sections of
$i^*(H^0(S,\mathcal{O}_S(k)\otimes I_x))$ also give an embedding of
$\overline C$ into ${\PP}^{r_k-1}$, since there is a natural injection
$
H^0(S,\mathcal{O}_S(k-2))\otimes H^0(S,\mathcal{O}_S(2)\otimes I_x)\subset
H^0(S,\mathcal{O}_S(k)\otimes I_x) ,
$
and sections of $H^0(S,\mathcal{O}_S(k-2))$ have no base points and fixed
components. Therefore, sections of $i^*(H^0(S,\mathcal{O}_S(k)\otimes I_x))$
separate points and tangent directions on $\overline C$.
$\hphantom{aaaa} \rule{6pt}{6pt}$

\vspace{0.1cm}

We say that a {\it linear subspace} $L_4$ of dimension $r-4$
{\it possesses a property} (L) {\it with respect to an embedding}
$X\subset {\PP}^r$, if the projection ${\pi}_{L_4}$
to ${\PP}^3$ from the centre $L_4$ maps $X$ onto a surface
$Y= {\pi}_{L_4}(X)$ with ordinary singularities.

As is known (see \cite{G-H} ), there is an open subset $U$ in the
Grassmannian $G_4=Gr(r-4,r)$, points of which correspond to linear subspaces
with the property (L).

\begin{pred}\label{P1.2}
If an embedding $X\subset {\PP}^r$ possesses the property (L), then there
exists a linear subspace $L_4$ with the property (L), which can be added to
a flag $L_1\supset L_3\supset L_4$ such that the pair $(L_1,L_3)$
possesses the property (L). In other words, there exists a projection to
${\PP}^3$, for which the embedding $Y\subset {\PP}^3$, where $Y$
is the image of $X$, possesses the property (L).
\end{pred}

\dvo Let $G_1= \check{\PP}^r $ be the dual space to ${\PP}^r$,
$G_3=G(r-3,r)$ be the Grassmann variety of linear subspaces $L_3$ of
dimension $r-3$, and
${\FF}={\FF}_{1,3,4} \subset \check{\PP}^r\times G_3\times G_4$ be the
variety of flags $L_1\supset L_3\supset L_4$. Let $X^{\vee}$ be the dual
variety, $W\subset X\times X^{\vee}\subset
{\PP}^r\times \check{\PP}^r$ be a closed subvariety
$W=\{ (x,L_1) : L_1\supset T_{X,x} \} $. Then the projection of
$W\to X^{\vee} $ is an isomorphism over
$X^{\vee}_0 =X^{\vee} \setminus Sing\; X^{\vee} ,$
$W_0 \simeq X^{\vee}_0 $.

Denote by $Z\subset X\times {\FF}$ a closed subvariety
$$
Z=\left\{ (x;L_1\supset L_3\supset L_4) \mid \; (x,L_1)\in W , L_3\ni x \right\} ,
$$
and by $Z_0\subset Z$ an open subset: $(x,L_1)\in W_0 $. Then $Z$ is an
irreducible variety. Consider a projection $Z_0 \to W_0$. The fibres are not
empty by the previous proposition, and each of the fibres contains an open
set of points $z$, for which the pair $(L_1,L_3)$ possesses the property (L)
(because the centres of projections for Lefschetz pencils form an open set).
Therefore, $Z$ contains an open set $Z_L$, for points of which the pair
$(L_1,L_3)$ possesses the Lefschetz property.

Obviously, the map $Z\to G_4$ is surjective. Therefore, $p_4^{-1}(U)$,
where $p_4$ is a projection of $Z$ to $G_4$, is a non empty Zariski open set.
Then $Z_L\cap p_4^{-1}(U)$ is not empty, and if
$(x;L_1\supset L_3\supset L_4)$ is a point of this set, then $L_4$
possesses  the desired property.
$\hphantom{aaa} \rule{6pt}{6pt}$

\vspace{0.1cm}

{\bf 1.7.} {\it Projecting in a neighbourhood of a generic point of the
double curve} $\Delta $.
Now let $Y\subset {\PP}^3$ has ordinary singularities along the double curve
$\Delta$ and isolated singularities $s_i\in Y\setminus\Delta ,
i=1,\ldots ,k$, which are double planes. Under the incidence varieties,
defined in the previous subsections, we mean the closures of the
corresponding varieties, initially defined for an open surface
$Y\setminus Sing\;Y $.

Consider $Y\times{\FF}$. In addition to notations in 1.3, let
$q_1$ and $q_2$ the projections of $Y\times G$ to $Y$ and $G$. Consider
the intersection $\tilde A =(\Delta\times G)\cap Z_3$.
Then the restriction of the projection $Y\times G \to Y$ to
$\tilde A,$ $\tilde A \to \Delta$, is a covering of degree 4 over a generic
point : at a point $y\in {\Delta}$ there are two asymptotic directions for
each of two branches of $Y$ at $y$. Therefore, $\tilde A$ is a curve. Set
$A= ({pr}'_{2})^{-1}(\tilde A)$. It is a ruled  surface. Set
${\Sigma}_{\Delta} =p_2(A)$. Then, if $\xi\notin{\Sigma}_{\Delta}$, we have
that for a generic point $y \in \Delta$ the lines
$l\ni\xi$ have at most simple contact with branches of $Y$.

Denote by ${\Sigma}_0$ the union of planes in ${\PP}^3$ composing the tangent
cones at the rest points of $\Delta $, including triple points and pinches,
and also at singular points $s_i\in Y$.

\vspace{0.1cm}

{\bf 1.8.} {\it Projecting in a neighbourhood of a triple point.}
If $\xi\notin{\Sigma}_0$, then in a neighbourhood of a point $y\in{\Sigma}_t$
all lines $l\ni\xi$ are transversal to each of the three branches of $Y$
at $y$, and therefore, locally these branches are mapped isomorphically onto
${\PP}^2$.

\vspace{0.1cm}

{\bf 1.9.} {\it Projecting in a neighbourhood of a pinch.}
In a neighbourhood of a pinch $y\in Y$ there are coordinates, by which
$Y$ is locally defined by an equation $u^2=vw^2$. The double curve
$\Delta\subset Y$ is a line $u=w=0$, and the tangent cone $C_{Y,y}$ to $Y$
at $y$ has an equation $u=0$.
In a neighbourhood of a pinch a normalization ${\CC}^2 \to Y$ is defined by
formulae $u=tw, v=t^2, w=w$. Since $X$ is non-singular and ${\pi}_L$ is a
finite map, we can assume that the projection ${\pi}_L$ is the normalization.
If a point $\xi$ does not belong to the tangent cone $C_{Y,y}$, then
the projection ${\pi}_{\xi}$ locally is a map of gedrr 2.
A projection $f: X \to {\PP}^2$ a neighbourhood of the preimage of a pinch
is a two-sheeted covering of non-singular varieties, and, hence, locally
is defined by equations $v=t^2, w=w$. Thus, the
curve $\bar R \subset Y$  goes through the pinch,
and pinches are projected to
non-singular points of the discriminant curve $B$.

\vspace{0.1cm}

{\bf 1.10.} {\it Normal forms of a generic projection at points of the
ramification curve.}

\begin{lem}[\cite{A}] \label{L3}
Let $(X,0)\subset (\CC^3,0)$ be a non-singular surface, and
$(\CC^3,0) \to (\CC^2,0)$ be a smooth morphism, the restriction of which
$f: X\to {\CC}^2$ is a finite covering of degree $\mu $. Then one can
choose local coordinates $x, y$ in $\CC^2$ and $x, y, z$ in $\CC^3$ such,
that $X$ is defined by an equation
$$
y=z^{\mu} + {\lambda}_1(x)z^{\mu -2} +\dots + {\lambda}_{\mu -2}(x)z ,
$$
and $f$ is a projection along $z$ axis.
\end{lem}

\dvo This is Lemma 1 in Arnol'd paper \cite{A}. It is obtained, if we
consider the covering $f$ as a 2-paremeter family of 0-dimensional
hypersurface singularities of multiplicity $\mu$, and, consequently, $f$
is induced by the miniversal deformation of the singularity of type
$A_{\mu -1}$.
\edvo

\vspace{0.1cm}
We proved that at a generic point $P$ of the ramification curve a projection
$f: X\to {\PP}^2$ is of degree $\mu =2$, and at isolated points it is of
degree $\mu =3$. By Lemma \ref{L3} for $\mu =2$ we obtain that at a generic
point of the ramification curve
a generic projection is equivalent to a projection of the surface $X: x=z^2$
to the $x,y$-plane, i.e. it is a fold. For $\mu =3$ we obtain

\begin{cor} \label{cor1}
For a non-singular surface $X$ a finite covering $X \to \CC^2$ of degree 3
locally is a projection to the $x,y$-plane of one of the surfaces
$$
y=z^3+x^kz  , \: k=1, 2,\dots , \; \mbox{or} \hphantom{a}
y=z^3 \; (k=\infty ) .
$$
In the case $k\ne \infty$ the ramification curve $C$ is reduced and has an
equation $3z^2+x^k=0$ in local coordinates $x, z$ on $X$. The curve $C$ is
non-singular only for $k=1$. The discriminant curve $B$ has an equation
$4x^{3k} +27y^2=0$, i.e. $B$ is a  cusp. It is an ordinary cusp only for
$k=1$.
\end{cor}
\dvo
By lemma \ref{L3} we have that $X$ is defined by an equation
$y=z^3 + {\lambda}_1(x)z$. We obtain the stated normal form of the covering
$f$, where $k$ is the order of vanishing of ${\lambda}_1(x)$ at the point
$x=0$. The ramification curve $C$ is defined by equation $J=0$, where
$J=3z^2 +x^k$ is the Jacobian of the covering $f$. The discriminant curve
$B$ is defined by 0th Fitting ideal $F_0\left( f_*{\cal O}_C\right)$
of the sheaf $f_*{\cal O}_C$. To obtain an equation of $B$ --- the generator
of the Fitting ideal, we need to take a finite presentation
$f_*{\cal O}_X \stackrel{J}{\to} f_*{\cal O}_X \to f_*{\cal O}_C \to 0$ of the sheaf
$f_*{\cal O}_C$, where
$\left( f_*{\cal O}_X\right)_0 =\CC\{ x,z\} =
\CC\{ x,y\}\cdot 1 \oplus \CC\{ x,y\}z \oplus \CC\{ x,y\}z^2$,
and to compute a determinant of the $\CC\{ x,y\}$-linear map $J$, which
is a multiplication by the Jacobian $J$.
\edvo

\vspace{0.1cm}
Now we show that for a generic projection the discriminant curve $B$ has
at most ordinary nodes and cusps.
Let $b\in B$ be a point corresponding to a bitagent $l$ under projecting
${\pi}_{\xi}: {\PP}^3\setminus \xi \to {\PP}^2 $ from a point $\xi $.
Let $l$ touches $Y$ at points $P_1$ and $P_2$, to which correspond branches
$B_1$ and $B_2$ of the discriminant curve $B$ at a point $b$. We have to
show that for a generic projection the point $b$ is a node, i.e. the branches
$B_1$ and $B_2$ have different tangents. Determine where does the centres
$\xi$ of "bad" projections lie. Let a line $\lambda \subset {\PP}^2$ is a
common tangent to branches $B_1$ and $B_2$ at a point $b$. Then the plane
${\pi}^{-1}_{\xi}(\lambda )$ is bitangent -- it touches the surface $Y$
at points $P_1$ and $P_2$. Consider the dual surface
$Y^{\vee}\subset \check{\PP}^3$. Then the point
$[{\pi}^{-1}_{\xi}(\lambda )]\in Sing\; Y^{\vee} = {\gamma}^{\vee}$.
Set $\gamma ={\tau}^{-1}({\gamma}^{\vee})$, where
$\tau : Y \to Y^{\vee}$ is the Gauss map. Let ${\Sigma}_u\subset
{\PP}^3$ be a ruled surface composed by lines $P_1P_2$, where
$P_1, P_2\in \gamma ,$ ${\tau}(P_1) ={\tau}(P_2) =
[{\pi}^{-1}_{\xi}(\lambda )]$. Then, if $\xi \notin {\Sigma}_u$, then
at points $b$, corresponding to bitangents $l$, the curve $B$ has at most
nodes.

Now let $b\in B$ be a point corresponding to a stationary tangent $l$
at a point $P\in Y$. As was noted above, in a neighbourhood of $P$ the
projection ${\pi}_{\xi}$ is equivalent to a projection of a surface
$y=z^3+x^kz $ to the $x,y$-plane. We have to show that for a generic
projection the exponent $k=1$. The fact is that, if $k>1$, then the point
$P$ is a planar point of the surface $Y$. Excepting the centres of
projection lying in tangent planes to
$Y$ at planar points, we obtain that in a neighbourhood of a point with
$\mu =3$ the projection $f$ is equivalent to a projection of a
surface $y=z^3+xz$ to the $x,y$-plane, i.e. it is a pleat.

\vspace{0.1cm}

{\bf 1.11.} {\it Projecting in a neighbourhood of an isolated double plane
singularity.}

\begin{lem} \label{L4}
If $(X,0) \subset (\CC^3,0)$ is an (isolated) double plane singularity
$z^2=h(x,y)$ , $ \pi : X\longrightarrow {\CC}^2$ be a projection from any
point $p\in {\CC}^3$, not lying in the tangent cone $z=0$, then the
ramification curve of $\pi$ is reduced, and the discriminant curve
$B\subset {\CC}^2$ is locally analytically isomorphic to the
singularity $h(x,y)=0$.
\end{lem}

\dvo The singularity $(X,0)$ is of multiplicity 2. Therefore, $\pi$ is a
covering of degree 2, and, consequently, is locally a projection of a
double plane $w^2=g(u,v)$ to the $(u,v)$-plane. Thus, the germs of
singularities $h(x,y)$ and $g(u,v)$ are stably isomorphic, and hence
isomorphic (\cite{AGV} , vol.1). $\hphantom{aa} \rule{6pt}{6pt}$

$\hphantom{aaaaaaaaaaaaaaaaaaaaaaaaaaaaaaaaaaaaaaaaaaaaaaaaaaaaaaaaaaaaaaa}
\rule{9pt}{9pt}$

\section{Local structure of fibre products.}

{\bf 2.1.} {\it Local structure of a generic covering.}
Let $f:X\to  {\PP}^2$ be a generic covering of the plane by a surface $X$
with A-D-E-singularities, and let $B\subset {\PP}^2$ be the discriminant
curve, $f^*(B) =2R+C$. Singular points $o \in {Sing}\,X $ will be called
{\it s-points } of the surface $X$ (from the word singularity). In a
neighbourhood of a s-point $o$ the covering $f$ is isomorphic to the projection
to $x,y$-plane of a surface $z^2 =h(x,y)$, where $h(x,y)$ has one of the
A-D-E-singularities. Singular points $o$ on $X$ correspond to singular points
of the same type as $o$ on $B$. With respect to $f$ non-singular points of
$X$ are partitioned into {\it r-points}
(from the word regular), at which the morphism $f$ is \'{e}tale, and
{\it p-points} (from the words singularity of projection) -- they are points
of the ramification curve $R$. A p-point is either a fold
(or a singular p-point of type $A_1$), in a neighbourhood of which
$f : (x,z) \longmapsto (x,y) ,\; y=z^2$,
or a pleat $o \in R\cap C$ (or a singular p-point of type $A_2$),
in a neighbourhood of which
$f : (x,z) \longmapsto (x,y) ,\; y=z^3 - 3xz$ (more details about this see
in section 2.4 below).

The singular points of $B$ \lz originated' from singular points ${Sing}\,X$
we call {\it s-points }. There are additional singular points of type
$A_1$ (nodes) and of type $A_2$ (cusps), which we call
{\it p-nodes } and {\it p-cusps}. Over a generic point $b\in B$ there lie
one fold and $N-2$ r-points; over p-node there lie two folds and $N-4$
r-points; over a p-cusp there lie one pleat and $N-3$
r-points; over a s-node or a s-cusp, as over any s-point $b\in B$, there
lie one singular point of $X$ and $N-2$ r-points.

{\bf 2.2.} {\it Types of points on the fibre product.} With respect to a
pair of generic coverings $f_1 :X_1 \to  {\PP}^2$ and
$f_2 :X_2 \to  {\PP}^2$ with the same discriminant curve
$B\subset {\PP}^2$ nodes and cusps on $B$ are partitioned by this time into
4 types:
{\it ss- ,} {\it sp- , ps-} and {\it pp-nodes } and {\it cusps.}
For example, a ps-node it is a node $b\in B$, such that there are two
folds on $X_1$ over $b$, and on $X_2$ over $b$ there is a singular
point of type $A_1$. The analogous terminology is used for the classification
of points $\bar x =(x_1,x_2) $ on the fibre product
$X^{\scriptscriptstyle \times} = X_1 \times_{{\PP}^2} X_2$ :
we say about rs-points, ss-points, sp-points, etc. For example, we say that
$\bar x$ is a ps-point of type $A_2$, if $x_1\in X_1$ is a
p-point of type $A_2$, and $x_2\in {Sing}\,X_2$ is a singular s-point
of type $A_2$.

In this section we describe the structure of a normalization
$\nu : X= \left( X^{\scriptscriptstyle \times }\right) ^{(\nu )} \to X^{\scriptscriptstyle \times} $
of the fibre product $X^{\scriptstyle \times}$. Denote by $g_1, g_2$ and
$f$ the morphisms of $X$ to $X_1, X_2$ and $\PP^2$.
Since the normalization
is defined locally, we can replace ${\PP}^2$ by a neighbourhood of the point
$0\in {\CC}^2$ and to assume that $X_1$ and $X_2$ are neighbourhoods of
points $x_1\in X_1$ and $x_2\in X_2$. We pass on to an item-by-item
examination of all possible types of points
$\bar x =(x_1,x_2) \in X^{\scriptscriptstyle \times} $. We do it up to the permutation of
factors $X_1$ and $X_2$.

At first we consider quite trivial cases.

\vspace{0.1cm}
\noindent 2.2.1. If $\bar x $ is a $r*$-point (where $*=r,s,p$), then
$X^{\scriptscriptstyle \times} $ at the point $\bar x$ is locally the same as $X_2$
at the point $x_2$, and $f^{\scriptscriptstyle \times} : X^{\scriptscriptstyle \times} \to
{\CC}^2 $ is locally the same as $f_2 : X_2 \to {\CC}^2$.

\vspace{0.1cm}
\noindent 2.2.2. If $\bar x $ is a {\it $2\times 2$-point}, i.e. $x_1$ and
$x_2$ are points of \lz double planes', $z_1^2=h(x,y) \; , z_2^2=h(x,y)$,
then $X^{\scriptscriptstyle \times} = X_1 \times_{{\CC}^2} X_2$ in a neighbourhood of the
point $\bar x =(x_1,x_2) $ is a surface in ${\CC}^4 \ni (x,y,z_1,z_2)$,
defined by equations $z_1^2=h(x,y) ,\; z_2^2=h(x,y)$. We obtain that
$z_1^2 =z_2^2$ and hence $X^{\scriptscriptstyle \times} =X_1^{\scriptscriptstyle \times} \cup
X_2^{\scriptscriptstyle \times}$, where
$X_1^{\scriptscriptstyle \times} : z_1^2 =h(x,y) , z_2=z_1 \; ,\;
X_2^{\scriptscriptstyle \times} : z_2^2 =h(x,y) ,\; z_1=-z_2$.  The surfaces
$X_1^{\scriptscriptstyle \times}$ and $X_2^{\scriptscriptstyle \times}$
meet along a curve
$z_1=z_2=0 ,\; h(x,y)=0$. We obtain that a normalization
$X=\widetilde {X^{\scriptscriptstyle \times}} =X_1^{\scriptscriptstyle \times} \sqcup
X_2^{\scriptscriptstyle \times} $ locally consists of two disjoint components
$X_1^{\scriptscriptstyle \times}$ and $X_2^{\scriptscriptstyle \times}$ isomorphic to
$X_1$ and $X_2$.

In particular, we obtain a description of the normalization in a
neighbourhood of a pp-point
$(x_1,x_2) \in X^{\scriptstyle \times}$ lying over a non-singular point of
$B$, $B: x=0$,

\vskip3cm 
\hskip2cm
\setlength{\unitlength}{0.1cm}
\begin{picture}(80,40)
  \epsfxsize=10.0cm
  \epsffile{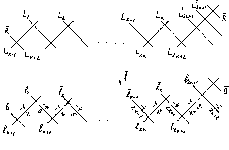}
\end{picture}

\begin{center}
Fig. 1 \\
\vspace{0.2cm}
$g_1^*(R_1) =R'+R'' , \; g_2^*(R_2) =R'+R'' $.
\end{center}

Every ss-point is a $2\times 2$-point. Thus, in a neighbourhood of a
ss-point the normalization has the same local structure as in the case of
a pp-point above: $X$ locally
consists of two disjoint components isomorphic to $X_1$ and $X_2$.

It remains to ezamine less trivial cases when $\bar x$ is a pp- or
sp-point of type $A_1$ or $A_2$. This is done in the following two
subsections.

{\bf 2.3.} {\it On fibre product of double planes.}
$\hphantom{aaaaaaaaaaaaaaaa}$

\noindent 2.3.1. The ordinary quadratic singularity -- the singularity  of type $A_1$
on a surface $X_0 : z^2=xy $ can be considered as a 2-sheeted covering of
the plane $f_0 : X_0 \to {\CC}^2 $ branched along a node $B: xy=0$.
As is known, the  singularity $X_0$ itself can be considered as a quotient
singularity under the action of cyclic group
${\ZZ}_2=\{ \pm 1\} , X_0 = X/{\ZZ}_2 $,
where $X={\CC}^2 \ni (z_1,z_2) $, and  a generator of ${\ZZ}_2$ acts by the
rule: $z_1 \longmapsto -z_1 , z_2 \longmapsto -z_2 $. The factorization
morphism $g_0 : X \to X_0 $ is defined by formulae
$$
x=z_1^2  ,\; y=z_2^2 , \; z=z_1z_2.
$$
We obtain a 4-sheeted covering $f=f_0\circ g_0 : X \to {\CC}^2 $,
which can be considered as a factorization under the action of the group
$G = {\ZZ}_2 \times {\ZZ}_2 $ on $X$. Then the factorization $g_0$
corresponds to a subgroup of order two ${\ZZ}_2 =G_0 = \{ (1,1),(-1,-1) \}$,
imbedded diagonally into $G$. In $G$ there are two more subgroups of order two:
$G_1 =\{ 1\} \times {\ZZ}_2$ and $G_2 = {\ZZ}_2 \times \{ 1 \} $.
Considering $X_1 ={\CC}_2 /G_1 \simeq {\CC}_2 $ and
$X_2 ={\CC}_2 /G_2 \simeq {\CC}_2 $, we obtain two more decompositions of
$f$ and a commutative diagram
\\

\unitlength=1.00mm
\special{em:linewidth 0.4pt}
\linethickness{0.4pt}
\begin{picture}(72.00,52.00)
\put(60.00,10.00){\makebox(0,0)[ct]{$(x,y)\in \mathbb C^2$}}
\put(80.00,50.00){\makebox(0,0)[cb]
{$X=\mathbb C^2 \ni (z_1,z_2) $}}
\put(55.00,30.00){\makebox(0,0)[rc]{$ (z_1,y)\in \mathbb C^2=X_1$}}
\put(78.00,30.00){\makebox(0,0)[lc]{$X_2= \mathbb C^2 \ni (x,z_2)\:$ ,}}
\put(66.00,48.00){\vector(0,-1){15.00}}
\put(71.00,40.00){\makebox(0,0)[rc]{$g_0 $}}
\put(66.00,28.00){\vector(0,-1){15.00}}
\put(71.00,22.00){\makebox(0,0)[rc]{$f_0 $}}
\put(64.00,48.00){\vector(-3,-4){10.67}}
\put(68.00,48.00){\vector(3,-4){10.67}}
\put(53.00,27.00){\vector(3,-4){10.67}}
\put(79.00,27.00){\vector(-3,-4){10.67}}
\put(57.00,20.00){\makebox(0,0)[rc]{$f_1 $}}
\put(75.00,20.00){\makebox(0,0)[lc]{$f_2 $}}
\put(58.00,42.00){\makebox(0,0)[rc]{$g_1$}}
\put(64.00,30.00){\makebox(0,0)[lc]{$X_0$}}
\put(74.50,42.00){\makebox(0,0)[lc]{$g_2$}}
\put(160.00,30.00){\makebox(0,0)[rc]{$(*_2)$}}
\end{picture}

\noindent where $g_1: y=z_2^2 , \: f_1: x=z_1^2 , \: g_2: x=z_1^2 , \:
f_2: y=z_2^2 $.

\noindent Denote by $B_1: x=0 ,\; B_2: y=0$ the branches of $B: xy=0 $, and
by $R': z_1=0 ,\; R'': z_2=0$  the branches of their proper transform
$z_1z_2 =0 $ on $X$.

The diagram $(*_2)$ shows that we can consider $X$ as a normalization in three
cases:

\vspace{0.1cm}
\noindent 2.3.2. $X$ is a normalization in a neighbourhood of a ps-point
of type $A_1$,
$\bar x \in X_1^{\scriptstyle \times} =X_1 {\times}_{{\CC}^2} X_0 $,

\vskip3cm 
\hskip3cm
\setlength{\unitlength}{0.1cm}
\begin{picture}(80,40)
  \epsfxsize=10.0cm
  \epsffile{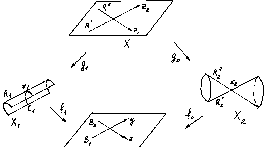}
\end{picture}

\begin{center}
Fig. 2
\end{center}

$$f_1^*(B_1) =2R_1 , \: f_1^*(B_2) =C_1 \; ; \: g_1^*(R_1) =R' , \:
g_1^*(C_1) =2R'' \; ; \:
$$
$$
f_0^*(B) =2R_2 =2(R'_2 +R''_2) \:  ; \;
g_0^*(R_2') =R' \; ; \: g_0^*(R_2'') =R'' \, . \:
$$
\noindent ($g_0$ is unramified outside the point $0\in X_0$ ).

\vspace{0.1cm}
\noindent 2.3.3. $X$ is a normalization in a neighbourhood of a
sp-point of type $A_1 ,$
$\bar x \in X_2^{\scriptstyle \times} =X_0 {\times}_{{\CC}^2} X_2 $
(the case symmetric to 2.3.2.)

\vspace{0.2cm}
\noindent 2.3.4. $X$ is a normalization in a neighbourhood of a
pp-point of $A_1 ,$
$\bar x \in X^{\scriptstyle \times} =X_1 {\times}_{{\CC}^2} X_2 $ ,
which is not a $2\times 2$-point,
$f_1: x=z_1^2 \: ,\; f_2: y=z_2^2 $.

\vskip2.5cm 
\hskip3cm
\setlength{\unitlength}{0.1cm}
\begin{picture}(80,40)
  \epsfxsize=10.0cm
  \epsffile{p2.eps}
\end{picture}

\vspace{0.3cm}
\begin{center}
Fig. 3
\end{center}

$$f_1^*(B) =2R_1+C_1 , \; g_1^*(R_1)=R' , \; g_1^*(C_1) =2R'' , \:
$$
$$
f_2^*(B) =2R_2+C_2 , \; g_2^*(R_2)=R'' , \; g_2^*(C_2) =2R'  .
$$

\medskip

Using 2.3.2-2.3.4, now we can describe a normalization $X$
over a neighbourhood of a node $b\in B $.

\vspace{0.1cm}
\noindent 2.3.5 Over a neighbourhood of a ps-node $b\in B \:$
(as well as a sp-node) a normalization of $X^{\scriptstyle \times}$
in a neighbourhood of a ps-point looks like as

\vskip2cm 
\hskip3cm
\setlength{\unitlength}{0.1cm}
\begin{picture}(80,40)
  \epsfxsize=10.0cm
  \epsffile{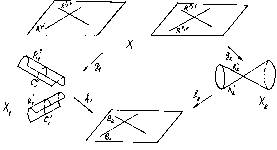}
\end{picture}

\begin{center}
Fig. 4
\end{center}
$$
g_1^*(R_1') =R^{\prime ,\prime} \; , \; g_1^*(C_1') =2R^{\prime ,\prime \prime}  \; , \;
g_1^*(R_1'') =R^{\prime \prime ,\prime} \; , \; g_1^*(C_1'') =2R^{\prime \prime ,\prime \prime}  \; , \;
$$
$$
g_2^*(R_2') =R^{\prime ,\prime}  +
R^{\prime \prime ,\prime  \prime} \; , \;
g_2^*(R_2') =R^{\prime ,\prime \prime}  +R^{\prime \prime ,\prime}  \;  \;
$$
\noindent On Fig. 4 the normalization in neighbourhoods of pr-, rs- and
rr-points of $X^{\scriptstyle \times} $ is not pictured.

\vspace{0.1cm}
\noindent 2.3.6. Over a neighbourhood of a pp-node $b\in B \:$
a normalization of $X^{\scriptstyle \times}$ in a neighbourhood of a
pp-point looks like as:

\vskip2cm 
\hskip2.8cm
\setlength{\unitlength}{0.1cm}
\begin{picture}(80,40)
  \epsfxsize=12.0cm
  \epsffile{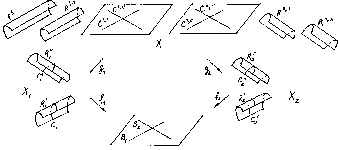}
\end{picture}

\begin{center}
Fig. 5
\end{center}
$$
g_1^*(R_1') =R^{\prime ,\prime}  + R^{\prime ,\prime \prime}  +
C^{\prime ,\prime}   \, , \;
g_1^*(R_1'') =R^{\prime \prime ,\prime}  +
R^{\prime \prime ,\prime \prime}  + C^{\prime \prime ,\prime}   \; , \;
$$
$$
g_2(R^{\prime ,\prime} ) =R_2' \; , \;
g_2(R^{\prime ,\prime \prime} ) =R_2' \; , \;
g_2(R^{\prime \prime ,\prime} ) =R_2'' \; , \;
g_2(R^{\prime \prime ,\prime \prime} ) =R_2'' \; , \;
$$
$$
g_2(C^{\prime ,\prime} ) =C_2'' \; , \;
g_2(C^{\prime \prime ,\prime} ) =C_2' \; , \;
\left( g_1(C^{\prime ,\prime \prime} ) =C_1' \; , \;
g_1(C^{\prime \prime ,\prime \prime} ) =C_1'' \right) \, . \;
$$

{\bf 2.4.} {\it On coverings of ${\CC}^2$ unbranched outside a cusp}
$B: y^2=x^3$.
To describe a normalization of the fibre product in a neighbourhood of a
sp- and pp-point of type $A_2$ in a natural context, we begin with reminding
of a small topic from singularity theory.

\noindent 2.4.1 {\it The singularity of cuspidal type of a map (pleat) and
the miniversal deformation of a singularity of type $A_2$.}
A cusp $(B,0) \subset ({\bf C}^2,0)$ is defined by a germ of
function $x^3-y^2$ stable equivalent to a germ of function $x^3$.
It is a simple singularity of type $A_2$.
It is interesting that a cusp (a singularity of type $A_2$ ) appears also
on the discriminant in the base of the miniversal deformation of the same
singularity of type $A_2$.

As is known, the miniversal unfolding of the function $t=z^3$ is
$$
{\CC}\times {\CC}^2 \to {\CC}\times {\CC}^2 \: , \:
(z,a_2,a_3) \longmapsto (z^3+a_2z+a_3, a_2, a_3) \, .
$$
The restriction of this map over $\{ 0\} \times {\CC}^2$ gives a miniversal
deformation $F$ of a zero-dimensional singularity $z^3=0$,
${\CC}^3\supset X \stackrel{F}{\longrightarrow} {\CC}^2 $. Here $X$
is a surface $z^3+a_2z+a_3 =0 $,  and $F$ is induced by projection onto
$(a_2,a_3)$. The surface $X$ is a graph of function $-a_3=z^3+a_2z $ ;
$z$ and $a_2$ are local coordinates on $X$,
\\

\unitlength=1.00mm
\special{em:linewidth 0.4pt}
\linethickness{0.4pt}
\begin{picture}(52.00,30.00)
\put(47.00,11.00){\makebox(0,0)[ct]{$(a_2,a_3)\in \mathbb C^2$}}
\put(45.00,30.00){\makebox(0,0)[rc]{$ (a_2,z)\in \mathbb C^2$}}
\put(47.00,31.00){\makebox(0,0)[lc]{$\stackrel{\sim}{\longrightarrow}$}}
\put(130.00,20.00){\makebox(0,0)[rc]{$,\, \, \, G:\left\{ \begin{array}{rcl}
a_2 & = & a_2 \\
-a_3 & = & z^3+a_2z.
\end{array}
\right. $}}
\put(56.00,27.00){\vector(0,-1){14.00}}
\put(60.00,20.00){\makebox(0,0)[rc]{$F $}}
\put(43.00,27.00){\vector(3,-4){10.67}}
\put(46.00,20.00){\makebox(0,0)[rc]{$G $}}
\put(54.00,30.00){\makebox(0,0)[lc]{$X\subset \mathbb C^3$}}
\end{picture}

We obtain a 3-sheeted covering $G: {\CC}^2 \to {\CC}^2 $,
the ramification curve of
which $R$ is defined by the equation $3z^2 +a_2 =0 $, and the discriminant
(branch) curve $B=G(R)$ is defined by equation
$$
4a_2^3 +27a_3^2 =0.
$$
To bring the equation of $B$ to the form $y^2=x^3$, we make a substitution
$$
a_2 =-3x \: , \: a_3=2y \, ,
$$
and denote ${\CC}^2 \simeq X$ by $X_3$, and $G$ by $f_3$.

\vskip2cm 
\hskip3cm
\setlength{\unitlength}{0.1cm}
\begin{picture}(80,40)
  \epsfxsize=10.0cm
  \epsffile{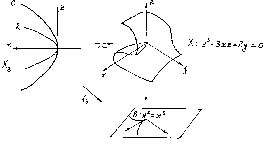}
\end{picture}

\nopagebreak
\begin{center}
Fig. 6
\end{center}

\noindent We obtain a 3-sheeted covering $f_3: X_3 \to {\CC}^2 $,
$$
f_3: \; x=x, \; y=-\frac12 (z^3-3xz).
$$
Then
$x^3-y^2 =x^3-\frac14(z^3-3xz)^2 =(x-z^2)^2(x-\frac14 z^2)
$
and, consequently,
$$
f_3^*(B) =2R+C  ,
$$
where $R: x=z^2$ is the ramification curve, and $C: x=\frac14 z^2 $. Note
that $C$ and $R$ are tangent of order two, $(C\cdot R)=2$.

By Lemma \ref{L3} the singular point of the covering $f_3$ is uniquely
characterized as a singular point of a 3-sheeted covering
$f: X \to \CC^2$ by a non-singular surface $X$, the discriminant curve of
which is an ordinary cusp.

\vspace{0.1cm}
\noindent 2.4.2 {\it The Vi\`{e}te map  $f_6$.}
We produce a well known regular covering of $\CC^2$ with group $S_3$ branched along a cusp
$B: y^2=x^3$, which appears to be a normalization of the fibre product
in a neighbourhood of a sp-point of type $A_2$. This covering naturally
appears in singularity theory.

Consider a quotient of the space ${\CC}^3$ under the action of permutation
group $S_3$. We get the Vi\`{e}te map
$$
v: {\CC}^3 \to {\CC}^3 \; \, ,\; (z_1,z_2,z_3) \longmapsto (a_1,a_2,a_3) ,
$$
where $(z-z_1)(z-z_2)(z-z_3) =z^3 +a_1z^2 +a_2z +a_3$ , i.e.
$$
a_1=-(z_1+z_2+z_3), \: a_2=z_1z_2+z_2z_3+z_3z_1, \; a_3=-z_1z_2z_3 .
$$
The map $v$ is a map of degree 6 unramified outside
$\Delta = \cup_{i\ne j} \{ z_i=z_j\} $,
and $v(\Delta )=D $ is defined by the discriminant of a polynomial of degree
three.

The action of $S_3$ on ${\CC}^3$ is reducible: ${\CC}^3$ is a direct sum
${\CC}^3 ={\CC}\oplus {\CC}^2$ of invariant subspaces -- of the line
${\CC} =\{ z_1=z_2=z_3\} $ and of the plane
${\CC}^2 =\{ z_1+z_2+z_3=0\} $. Consider the restriction of $v$ to this plane
${\CC}^2$,
$$
(z_1,z_2,z_3) \in {\CC}^3  \supset
\{z_1+z_2+z_3=0 \} = {\CC}^2 \longrightarrow {\CC}^2 = \{ a_1=0\}
\subset {\CC}^3 \ni (a_1,a_2,a_3) \, .
$$
\noindent Set ${\CC}^2 \cap \Delta =L$ , ${\CC}^2 \cap D =B$.
Then $L$ consists of three lines
$$
L=L_1+L_2+L_3 , \; \mbox{where} \: L_i : z_j=z_k ,\; z_1+z_2+z_3 =0 \; , \;
\{ i,j,k\} =\{ 1,2,3\} \, ,
$$
\noindent and the curve $B: 4a_2^3 +27a_3^2 =0 $ is defined by the
discriminant of the polynomial $z^3+a_2z+a_3 $. Since
${\pi}_1({\CC}^2 \setminus L ) ={\pi}_1({\CC}^3 \setminus \Delta ) \; , \;
{\pi}_1({\CC}^2 \setminus B ) ={\pi}_1({\CC}^3 \setminus D )= {Br}_3 $,
we  obtain
a covering $v: {\CC}^2 \to {\CC}^2 $ of degree 6 unbranched apart from $B$.
Denote this map by $f_6$. In coordinates $x, y $, where
$a_2=-3x$, $a_3 =2y$, this map
$$
{\CC}^2 =
\{z_1+z_2+z_3=0 \} = X_6
\stackrel{f_6}{\longrightarrow} {\CC}^2 \ni (x,y)
$$
is defined by formulae
$$
f_6 : \: x=-\frac13 (z_1z_2+z_2z_3+z_3z_1) \; , \; y =-\frac12 z_1z_2z_3 ,
$$
the discriminant $B$ has equation $y^2=x^3$, and $f^*(B) =2L=
2L_1+2L_2+2L_3 $ (it is easy to see that
$x^3-y^2 =\frac{1}{4\cdot 27} (z_2-z_1)^2(z_3-z_2)^2(z_1-z_3)^2$
under condition $z_1+z_2+z_3 =0$ ).

Consider a two-sheeted covering unbranched outside $B$
$$
(x,y,w)\in {\CC}^3 \supset X_2
\stackrel{f_2}{\longrightarrow} {\CC}^2 \ni (x,y) \; ,
$$
where $X_2$ is defined by equation $w^2=x^3-y^2$, and $f_2$ is induced by
projection. Such a structure has a generic covering $f: X \to {\PP}^2$
in a neighbourhood of a s-point of type $A_2$.

\begin{lem}[\cite{C}] \label{L5}
If $f: (X,0) \to (\CC^2,0)$ is a finite covering by a normal irreducible
surface $X$, unbranched outside an ordinary cusp $B\subset \CC^2$, and the
ramification curve of which is reduced, i.e. $f^*(B) =2R+C$ ($R$ and $C$
reduced curves), then $f$ is equivalent to one of the coverings $f_2, f_3$
and $f_6$.
\edvo
\end{lem}
The proof is obtained by means of studying the possible monodromy
homomorphisms $\rho : {\pi}_1 \to S_N$, where ${\pi}_1 ={\pi}_1(\CC^2
\setminus B) = Br$ is the fundamental group of a cusp, and $N=\deg f$.

We obtain one more characterization of the covering $f_3$ as a finite
covering $f: (X,0) \to (\CC^2,0)$ by a normal irreducible surface, unbranched
outside a cusp $B$,  and with a reduced and non-singular ramification
curve $R$.

\vspace{0.1cm}
\noindent 2.4.3 {\it Description of a normalization of the fibre product
in a neighbourhood of a sp-point of type $A_2$.}
The map $f_6$ factors through the maps $f_2$ and $f_3$, and we have a
commutative diagram
\\

\unitlength=1.00mm
\special{em:linewidth 0.4pt}
\linethickness{0.4pt}
\begin{picture}(72.00,52.00)
\put(72.00,10.00){\makebox(0,0)[ct]{$(x,y)\in \mathbb C^2\supset B : y^2=x^3$}}
\put(86.00,50.00){\makebox(0,0)[cb]
{$X_6=\{z_1+z_2+z_3=0 \} $}}
\put(55.00,30.00){\makebox(0,0)[rc]{$\{ w^2=x^3-y^2 \} =X_2$}}
\put(78.00,30.00){\makebox(0,0)[lc]{$X_3= \{z^3-3xz^2+2y=0 \} \, , $}}
\put(66.00,48.00){\vector(0,-1){35.00}}
\put(64.00,48.00){\vector(-3,-4){10.67}}
\put(68.00,48.00){\vector(3,-4){10.67}}
\put(53.00,27.00){\vector(3,-4){10.67}}
\put(79.00,27.00){\vector(-3,-4){10.67}}
\put(57.00,20.00){\makebox(0,0)[rc]{$f_2 $}}
\put(75.00,20.00){\makebox(0,0)[lc]{$f_3 $}}
\put(58.00,42.00){\makebox(0,0)[rc]{$g_3$}}
\put(67.00,30.00){\makebox(0,0)[lc]{$f_6$}}
\put(74.50,42.00){\makebox(0,0)[lc]{$g_2$}}
\put(160.00,30.00){\makebox(0,0)[rc]{$(*_3)$}}
\end{picture}
\newline

\noindent where $g_2$ and $g_3$ are defined by formulae:
$x$ and $y$ are defined by the same formulae as $f_6$, and $z=z_1$ for
$g_2$, and $w=\frac{1}{6\sqrt 3}
(z_2-z_1)(z_3-z_2)(z_1-z_3)$ for $g_3$.
It is easy to see that $g_3$ is a factorization under
the action of a cyclic group ${\ZZ}_3 ={\cal A}_3 \subset S_3$,
$\; X_2 =X_6 / {\cal A}_3$, and $g_2$ is a factorization under
the action of a cyclic group of order two
${\ZZ}_2 \simeq S_2 =\{ (1),(2,3)\} \subset S_3 $.

By the property of universallity of fibre products we have a morphism
$X_6 \to X_2 {\times}_{{\CC}^2} X_3$. The fibre product
$X_2 {\times}_{{\CC}^2} X_3$ is irreducible, since each its component $Z$ is
mapped onto $X_2$ and $X_3$, and, therefore, the degree of $Z\to \CC^2$
have to be divided by 2 and 3, i.e. have to be equal to 6. Thus, $X_6$ is a
normalization of $X_2 {\times}_{{\CC}^2} X_3$, and the diagram $(*_3)$
describes a normalization of the fibre product in a neighbourhood
of a sp-point of type $A_2$.
\newpage 

The diagram $(*_3)$ can be visually-schematic presented as follows

\vskip2.5cm 
\hskip3cm
\setlength{\unitlength}{0.1cm}
\begin{picture}(20,40)
  \epsfxsize=10.0cm
  \epsffile{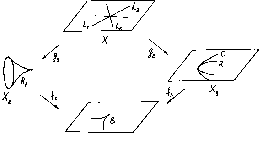}
\end{picture}
\begin{center}
Fig. 7
\end{center}

Direct computations show that $x-z^2=\frac13 (z_2-z_1)(z_1-z_3)$ ,
and $x-\frac14 z^2=\frac{1}{12} (z_3-z_2)^2$, i.e.
$$
g^*_2(R) =L_2 +L_3 \; , \; g^*_2(C) =2L_1 .
$$
\noindent And, besides, $g_3^*(R_1)=L_1+L_2+L_3$.

\vspace{0.1cm}

\noindent 2.4.4 {\it Description of a normalization of the fibre product
in a neighbourhood of a pp-point of type $A_2$.}
Let $x_1 \in X_1$ and $x_2\in X_2$ be p-points of type $A_2$
for $f_1$ and $f_2 $,
$
f^*_1(B) =2R_1+C_1 \: , \: f^*_2(B) =2R_2+C_2
$.
In this case the 3-sheeted coverings $f_1$ and $f_2$ are the same
(equivalent), and the monodromy homomorphisms
${\varphi}_1 ,{\varphi}_2: {\pi}_1 ={\pi}_1(\CC^2 \setminus B, y_0) \to S_3$
are epimorphic. The fibre $(f^{\scriptscriptstyle \times})^{-1}(y_0)$ of the 9-sheeted covering
$f^{\scriptscriptstyle \times}:
X^{\scriptscriptstyle \times} = X_1 {\times}_{{\CC}^2} X_2 \to \CC^2$
consists of pairs
$f_1^{-1}(y_0)\times f_2^{-1}(y_0) =\{ (i,j) , \: 1\le i, j \le 3 \}$,
and the monodromy homomorphism is (equivalent to) a diagonal homomorphism
$\varphi : {\pi}_1 \to S_3\times S_3\subset S_9$. Since ${\varphi}_i$ are
epimorphic, the fibre $(f^{\scriptscriptstyle \times})^{-1}(y_0)$ consists
of two orbits w.r.t.
the action of ${\pi}_1$ --- the orbit of the point $(1,1)$, which consists
of 3 elements, and the orbit of the point $(1,2)$, which consists of 6 elements.
From this and from Lemma \ref{L5} it follows that in a neighbourhood of the
$\bar x =(x_1,x_2)$ a normalization $X$ of the product
$X^{\scriptscriptstyle \times}$ consists of two disjoint components
$X=X_3\coprod X_6 $, and
on $X_3$ the morphism $f: X\to \CC^2$
coincides with $f_3$, the morphisms  $g_1$ and $g_2$ are isomorphisms, and
on $X_6$ the morphism $f=f_6$,  the morphisms $g_1$ and $g_2$ are the same
as $g_2$ in the diagram $(*_3)$
\newpage 
\vskip2cm 
\hskip3cm
\setlength{\unitlength}{0.1cm}
\begin{picture}(80,40)
  \epsfxsize=10.0cm
  \epsffile{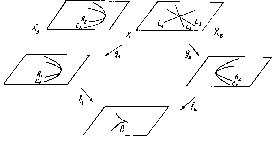}
\end{picture}

\begin{center}
Fig. 8
\end{center}

\noindent There are 4 curves on
$X^{\scriptscriptstyle \times}$:
 $C_1 {\times}_{B} C_2$, $R_1 {\times}_{B} C_2$,
$C_1 {\times}_{B} R_2$, $R_1 {\times}_{B} R_2$, preimages of which on the
normalization $X$ are $C_3$, $L_2$, $L_1$, and $L_3$, $R_3$. Under such a
numeration of the lines $L_i$ we have
$$
g^*_1(R_1) =R_3 + L_2 +L_3 \; , \; g^*_1(C_1) =C_3 +2L_1 ,
$$
$$
g^*_2(R_2) =R_3 + L_1 +L_3 \; , \; g^*_2(C_2) =C_3 +2L_2 .
$$

\vspace{0.1cm}

\noindent 2.4.5 {\it A lift of the diagram $(*_3).$}
Consider the diagram $(*_3) .$ For computation of intersection numbers in
$\S 5$ we need to resolve the singular point of type $A_2$ on the surface
$X_2$, and to \lz disjoint' the curves $L_2$ and $L_3$ on $X$. A resolution
of the singular point of type $A_2 $,  as of any \lz double plane', can be
obtained, if we firstly take an imbedded resolution
$\sigma : \bar {{\CC}^2} \to {\CC}^2$ of the branch curve
$B \subset {\CC}^2 $, and then take a normalization of
$X_2 {\times}_{{\CC}^2}{\bar {\CC}^2}$.

Actually we'll make more -- we lift the whole of the diagram $(*_3)$ on
$\bar {{\CC}^2}$.

1) The singular point of $B$ is resolved by one $\sigma$-process ${\sigma}_1$.
It is enough for the resolving of the singular point on $X_2$, but to resolve
the total transform of $B$ up to a divisor with normal crossings, one need
two more $\sigma$-processes. We picture the resolution process schematically
by \lz drawing' the total transform of the curve $B$. Denote by $E_i$ the curve
glued in under the i-th $\sigma$-process, and also its proper transform under
subsequent $\sigma$-processes.

\vskip1.2cm 
\hskip2.5cm
\setlength{\unitlength}{0.1cm}
\begin{picture}(80,40)
  \epsfxsize=12.0cm
  \epsffile{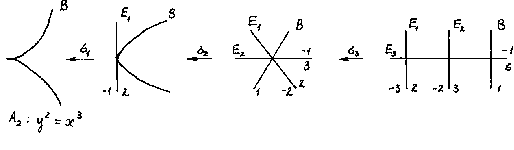}
\end{picture}

\begin{center}
Fig. 9
\end{center}

\noindent Along each curve we indicate two numbers: the negative is the
self-intersection number, the positive is its multiplicity in the total
transform of the curve $B$.

2) Denote by $\sigma : \bar {{\CC}^2} \to {\CC}^2$ the composition
${\sigma}_3 \circ {\sigma}_2 \circ {\sigma}_1 $. We add on the diagram
$(*_3)$ over $\bar {{\CC}^2}$ and obtain a diagram as follows, on which
all morphisms on the right face are finite coverings.

\renewcommand{\theequation}{\fnsymbol{equation}}

\begin{equation}
\xymatrix
{
 & X_6 \ar[ld]_{g_3} \ar@{-->}[rdd]^(.35){g_2} & & \overline{X}_6 \ar[ll]_{\sigma_6}
\ar[ld]_{\overline g_3} \ar[rdd]^{\overline g_2} &   \\
X_2 \ar[rdd]_{f_2} & & \overline X_2 \ar[ll]^{\sigma_2}
\ar[rdd]^(.35){\overline f_2} & & \\
& & X_3 \ar@{-->}[ld]_{f_3} & & \overline X_3 \ar@{-->}[ll]_{\sigma _3}
\ar[ld]^{\overline f_3} \\
& \mathbb C^2 & & \overline{\mathbb C}^2 \ar[ll]_{\sigma } &
}
\end{equation}

The right square of the diagram (\theequation ) is obtained as a fibre
product $(*_3) {\times}_{\CC^2}\bar\CC^2$, i.e. $\bar X_i$ are normalizations
of $X_i {\times}_{\CC^2} \bar\CC^2$, and morphisms are induced by morphisms
of the diagram $(*_3)$ and projections. We describe how one can construct the
diagram (\theequation ) not uniformly as a normalization of the lift, but
step-by-step.
To facilitate the following of the description we begin with the final
picture. We draw the right square of the diagram
(\theequation ) by replacing
the varieties at its vertices by the total transforms of the curve $B$

\vskip3cm 
\hskip2cm
\setlength{\unitlength}{0.1cm}
\begin{picture}(100,80)
  \epsfxsize=12.0cm
  \epsffile{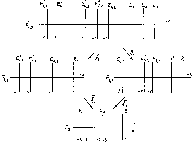}
\end{picture}

\begin{center}
Fig. 10
\end{center}

\noindent The rule of notation is as follows. The exceptional curves
$E_1, E_2, E_3 $ on $\bar {{\CC}^2}$ are already denoted. Under double
indexing $E_{i,j}$ the first index indicates the variety $X_i$, where
$E_{i,j}$ lies, and the second index indicates to what curve $E_j$
the curve $E_{i,j}$  is mapped on $\bar {{\CC}^2}$.

3) We begin a description of the diagram (\theequation )
with $\bar {X_6}$
(\lz from the top').
To disjoint the lines $L_i$, we make $\sigma$-process with centre at the
point $0\in X_6 ={\CC}^2 = \{ z_1+z_2+z_3 =0\} \subset {\CC}^3$.
By this the curve $E_{6,3} ={\PP}^1 =\{ t_1+t_2+t_3 =0\}
\subset {\PP}^2 \ni (t_1:t_2:t_3)$ is glued, and we obtain a variety $X_6'$.
The action of $S_3$ on $X_6$ is extended to $X_6'$ and, in particular, to
${\PP}^1$. On ${\PP}^1$ there are 8 exceptional points forming exceptional
orbits:
$$
p_1 =E_{6,3} \cap L_1 =(-2:1:1) , \;
p_2 =E_{6,3} \cap L_2 =(1:-2:1) , \;
p_3 =E_{6,3} \cap L_3 =(1:1:-2) ; \;
$$
$$
P_1=(0:1:-1) , \; P_2=(1:0:-1) , \; P_3=(1:-1:0) ; \;
$$
$$
Q_1=(1:\zeta :{\zeta}^2) , \; Q_2=(1:\bar{\zeta} :\bar{{\zeta}^2}) , \;
$$
where $ \zeta =\sqrt[3]{1}$ is a primitive root, and $\bar {\zeta} ={\zeta}^2$.
Denote by $\xi =(123) $ a generator of the cyclic group of order three
${\ZZ}_3 ={\cal A}_3 =\{ (1), (123), (132) \} \subset S_3 $, and by
$\varepsilon =(23) $ a generator of the cyclic group of order two
${\ZZ}_2 =S_2 =\{ (1), (23) \} \subset S_3 $. Then
$$
\xi (p_1) = p_2 , \; \xi (p_2) = p_3 , \; \xi (p_3) = p_1 ; \quad
\xi (P_1) = P_2 , \; \xi (P_2) = P_3 , \; \xi (P_3) = P_1 ;
$$
$$
\xi (Q_1) = ({\zeta}^2:1:\zeta) =(1:\zeta :{\zeta}^2) =Q_1 , \;
\xi (Q_2) = ({\zeta}:1:{\zeta}^2) =(1 :{\zeta}^2:\zeta) =Q_2 ;
$$
$$
\varepsilon (p_1) =p_1 , \; \varepsilon (p_2) =p_3 , \;
\varepsilon (p_3) =p_2 ; \;
\varepsilon (P_1) =P_1 , \; \varepsilon (P_2) =P_3 , \;
\varepsilon (P_3) =P_2 ; \; \varepsilon (Q_1)=Q_2 .
$$
If we take a quotient $X_6'$ under the action of ${\bf Z}_3 ={\cal A}_3$,
then the stationary points $Q_1$ and $Q_2$ give two quotient singularities
on $X_2' =X_6' / {\cal A}_3$, resolving of which $X_2'' \to X_2' $
glues the curves $E_{2,1}'$ and $E_{2,1}''$  with
$\left( {E_{2,1}'}^2 \right) =-3 $, $\left( {E_{2,1}''}^2 \right) =-3 $.
To lift $X_6' \to X_2' $ onto $X_2''$, we have to blow up the points $Q_1$
and $Q_2$ , $X_6'' \to X_6' $, and by this we obtain $X_6'' \to X_2''$,

\vskip0.5cm 
\hskip2.5cm
\setlength{\unitlength}{0.1cm}
\begin{picture}(80,40)
  \epsfxsize=12.0cm
  \epsffile{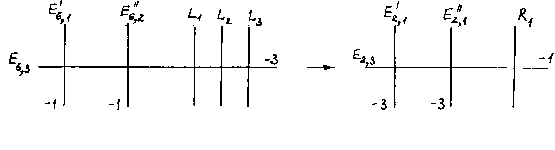}
\end{picture}

\begin{center}
Fig. 11
\end{center}

4) The map $f_2$ is a factorization under the cyclic group
${\ZZ}_2 =S_2$. The action extends to $X_2''$. The stationary point on
$E_{2,3}$ -- the image of the point $P_1$ on $E_{6,3}$ gives a singular
point of type  $A_2$ on $X_2'' / {\ZZ}_2 \:$
$ \left( =\bar {{\bf C}^2}' \right) $. A resolution of this point glues a
(-2)-curve $E_2$, and we obtain $\bar {\CC}^2 $. To lift
$X_2'' \to X_2'' / {\ZZ}_2 $ onto the resolution
$\bar {\CC}^2$, we have to blow up a point on $X_2''$. By this a
(-1)-curve $E_{2,2}$ is glued, and we obtain $\bar X_2 $. To obtain
$\bar g_3: \bar X_6 \to \bar X_2 $, we have to perform 3 $\sigma$-processes
with centres at points $P_1$, $P_2$, $P_3$ on $X_6'' $, by which three
lines $E_{6,2}'$, $E_{6,2}''$ and $E_{6,2}'''$ are glued. We obtain the left
side $\bar g_3$ and $\bar f_2$ of the right square of
diagram (\theequation ), pictured on
Fig. 10. Note thay the map $\bar g_3$ is ramified along the curves
$E_{6,1}'$ and $E_{6,1}''$, and the map $\bar f_2$ is ramified along the
curves $E_{2,2}$ and $R_1$.

We can blow down the (-1)-curve $E_{2,2}$ on $\bar X_2$, and then to blow
down the (-1)-curve $E_{2,3}$.  By this we obtain a minimal resolution of
the singular point of type $A_2$ on $X_2$,

\vskip0.5cm 
\hskip2.4cm
\setlength{\unitlength}{0.1cm}
\begin{picture}(80,40)
  \epsfxsize=10.0cm
  \epsffile{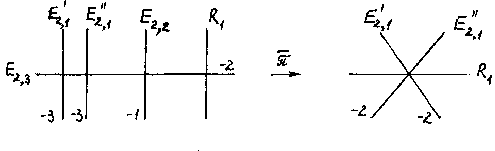}
\end{picture}

\begin{center}
Fig. 12
\end{center}

5) The map $\bar g_2 $ is a factorization under the group
${\bf Z}_2 =S_2 = \{ (1),(23)\} $. We obtain the surface
$\bar X_3 =\bar X_6 /S_2 , \; \bar g_2 : \bar X_6 \to \bar X_3 $.
The map $\bar g_2$ is ramified along the curves $E_{6,2}'$ and $L_2$,
which are mapped onto $E_{3,2}'$ and $C$ correspondingly.
The diagram is completed by the map  $\bar f_3: \bar X_3 \to \bar {\CC}^2 $.
The surface $\bar X_3$ is obtained from $X_3 ={\CC}^2 $, if we at first
blow up the point of tangency of curves $C$ and $R$ gluing $E_{3,2}'$;
then we blow up the point of intersection of $C$ and $R$ gluing $E_{3,3}$;
finally, we blow up two more points on $E_{3,3}$:

\vskip0.3cm 
\hskip2.3cm
\setlength{\unitlength}{0.1cm}
\begin{picture}(80,40)
  \epsfxsize=12.0cm
  \epsffile{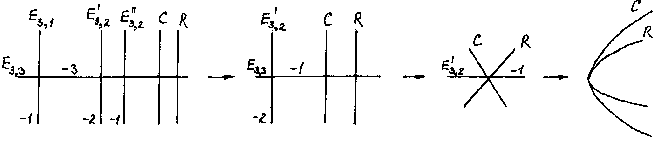}
\end{picture}

\begin{center}
Fig. 13
\end{center}

\section{The canonical cycle of a Du Val singularity}

We intend to apply Hodge index theorem to obtain the basic inequality for
generic coverings of ${\PP}^2$ by surfaces with A-D-E-singularities.
For this we need intersection theory and, therefore, a resolution of
singularities of $X$. In this section we examine the local situation and
find out how the resolution affects
the canonical class and the ramification curve.

\vspace{0.1cm}

{\bf 3.1.} {\it Definition of canonical cycle.}
Let $(X,x)$ be
a 2-dimensional A-D-E-singularity. Let $\pi : \bar X \to X $ be a
minimal resolution, $L={\pi}^{-1}(x)$ be the exceptional curve. As is known,
the canonical class $K_{\bar X}$ is trivial in a neighbourhood of $L$,
that is we can choose a divisor in $K_{\bar X}$ with a support not
intersecting $L$. In other words, there is a differetial form $\omega$
on $\bar X$, which has neither poles nor zeroes in a neighbourhood of $L$.
Such a form can be obtained, for example, as follows. As is known, $(X,x)$
is a quotient singularity, $X= {\CC}^2 / G$, where $G\subset SL(2,{\CC}) $ .
The form $du\wedge dv $ on ${\CC}^2 \ni (u,v) $ is invariant w.r.t. $G$
and it defines a form on $X \:$  (${\varphi}^*(\omega ) = du\wedge dv $,
where $\varphi : {\CC}^2 \to X $). Hence, the divisor
$(\omega ) =\sum k_iL_i $. Since $L_i$ are (-2)-curves,
$\bigl( L_i\cdot (\omega )\bigr) =0$,  and we obtain $(\omega ) =0 $.

  On the other hand, $(X,x)$ can be considered as a double plane,
that is as a 2-sheeted covering $X\stackrel{f}{\longrightarrow} Y$
of the plane $Y={\CC}^2$ (locally). Let $z^2=h(x,y)$ be an equation of
$(X,x)$, $B: h(x,y)=0$ be the discriminant curve, $f^{-1}(B)=R$,
defined by the equation $z=0 $,
be the ramification curve. We can consider the differential form
$\omega =f^*(dx\wedge dy) $ lifted from $Y$. Then on $\bar X$ the divisor
$(\omega ) =(z) =\bar R + Z $, where $\bar R \subset \bar X$ is the proper
transform of $R$ , $Z =\sum {\gamma }_iL_i $ is a cycle on
$L={\pi}^{-1}(x) $. We shall say that $Z$ is the {\it canonical cycle}
of a 2-dimensional A-D-E-singularity. Thus, $-Z$ is a cycle on the exceptional curve
$L$, which is equivalent to the ramification curve $\bar R$
in a neighbourhood of $L$. Let us calculate the canonical cycle for all
A-D-E-singularities.

\vspace{0.2cm}

{\bf 3.2.} {\it On resolution of double planes.}
As for any double plane, a resolution of an A-D-E-singularity can be
obtained by means of a resolution of the discriminant curve
$B \subset Y={\CC}^2 $, $B: h(x,y)=0$. Let $\sigma : \bar Y \to Y $ be a
composition of $\sigma$-processes, such that the total transform of $B$
is a divisor with normal crossings. Let
${\sigma}^*(B) = \bar B +\sum_{i=1}^{r} {\alpha}_il_i $, where $\bar B$
is the proper transform of $B$,  $l_i \simeq {\bf P}^1 , \; i=1,\ldots ,r $,
are the exceptional curves, as well as their proper transforms,
glued by $\sigma$-processes. Let $\bar X$ be the normalization of
${\bar Y}\times_YX$, and $\bar f$ and $\pi$ be induced by projections,

\unitlength=1mm
\special{em:linewidth 0.4pt}
\linethickness{0.4pt}
\begin{picture}(93.00,28.00)(10.00,4.00)
\put(89.50,28.00){\makebox(0,0)[rc]{$\pi ^{-1}(x) =
L=L_1\cup \ldots \cup L_r\subset \bar X$}}
\put(93.00,28.00){\vector(1,0){15.00}}
\put(110.00,28.00){\makebox(0,0)[lc]{$X\supset R \ni x , \; R: z=0$}}
\put(100.00,29.00){\makebox(0,0)[cb]{$\pi $}}
\put(87.50,25.00){\vector(0,-1){12.00}}
\put(89.00,10.00){\makebox(0,0)[rc]{$\sigma ^*(B)=\overline B+
\sum_{i=1}^{r}\alpha _il_i\subset\bar Y$}}
\put(93.00,10.00){\vector(1,0){15.00}}
\put(110.00,10.00){\makebox(0,0)[lc]{$Y \supset B .$}}
\put(112.00,25.00){\vector(0,-1){12.00}}
\put(113.00,19.00){\makebox(0,0)[lc]{$f$}}
\put(86.50,19.00){\makebox(0,0)[rc]{$\bar f $}}
\put(100.00,9.00){\makebox(0,0)[ct]{$\sigma$}}
\put(167.00,19.00){\makebox(0,0)[ct]{$(\sharp )$}}
\end{picture}
\newline

\noindent Set ${\bar f}^{-1}(l_i) = L_i$. The curve $L_i$ is either
irreducible or consists of two components $L_i =L_i' + L_i'' $, where
$L_i' \simeq {\PP}^1 $ , $L_i'' \simeq {\PP}^1 $. The mapping $\bar f$ is a
2-sheeted covering branched along the curve
$\bar B +\sum_{{\alpha}_i - \mbox{odd}} l_i $.
To be more graphic we denote the curves $l_i$, for which ${\alpha}_i$ are
odd, also by $\bar l_i $, and $L_i$ -- respectively by $\bar L_i$.
The surface $\bar X$ has singularities of type $A_1$ over nodes of the
branch curve $\bar B +\sum \bar l_i $. If this curve is non-singular,
that is, a disconnected union of components (one can reach this by performing
one additional $\sigma$-processes for each node), then $\bar X$  is
non-singular and is a resolution of the singularity $(X,x)$.
Let $\bar R $ be the proper transform of $R$ w.r.t. $\pi \:$
(= the proper transform of $\bar B$ w.r.t. $\bar f$ ). We have
${\bar f}^*({\bar l_i}) =2{\bar L}_i $, if ${\alpha}_i$  is odd, and
${\bar f}^*({l_i}) ={L}_i $, if ${\alpha}_i$  is even. We have
$$
\bigl( (\sigma \circ \bar f)^* h(x,y) \bigr) = (z^2) =
2\bar R +\sum_{{\alpha}_i - \mbox{odd}} 2{\alpha}_i \bar L_i +
\sum_{{\alpha}_i - \mbox{even}} {\alpha}_i L_i
$$
and, consequently, $(z) =\bar R +Z $, where
$$
Z= \sum_{{\alpha}_i - \mbox{odd}} {\alpha}_i \bar L_i +
\sum_{{\alpha}_i - \mbox{even}} \frac12 {\alpha}_i L_i \; .
$$
Let us compute the cycle $Z$ for each type of A-D-E-singularities
(despite of abundance of papers concerning  Du Val singularities, the
authors do not know any of them, where the cycle $Z$ is
written out; so we have to perform these computations).

\vspace{0.1cm}

{\bf 3.3.} {\it Computation of the canonical cycle.}
Consider the minimal resolution of each type of A-D-E-singularities
described above. The following lemma contains the results of computations of
${\sigma}^*(B)$, of the exceptional curve ${\pi}^{-1}(x)=L$ and of the
canonical cycle $Z$.

\begin{lem}
Below we picture schematically the total transform
${\sigma}^*(B) = \bar B +\sum_{i=1}^{r} {\alpha}_il_i \:$ (near each curve
$l_i$ a positive number ${\alpha}_i$ and a negative number
$(l_i^2)$ are shown), and over it we picture the curve ${\pi}^{-1}(R)$,
consisting of $\bar R$ and (-2)-curves, and besides we write down the
canonical cycle $Z$:

\noindent 1) The singularity  $A_{2k-1}: y^2=x^{2k} ,\; k\ge 1$,

\vskip3cm 
\hskip2.5cm
\setlength{\unitlength}{0.1cm}
\begin{picture}(120,60)
  \epsfxsize=10.0cm
  \epsffile{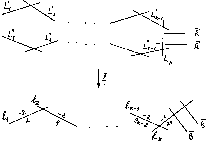}
\end{picture}

\begin{center}
$Z=L_1 +2L_2+\ldots +kL_k \, ;$
\end{center}

\newpage \noindent 2) The singularity  $A_{2k}: y^2=x^{2k+1} ,\; k\ge 1$,

\vskip3cm 
\hskip2.5cm
\setlength{\unitlength}{0.1cm}
\begin{picture}(120,60)
  \epsfxsize=10.0cm
  \epsffile{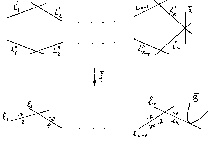}
\end{picture}

\begin{center}
$Z=L_1 +2L_2+\ldots +kL_k \, ;$
\end{center}

\noindent 3) The singularity   $D_{2k+2}: x(y^2+x^{2k}) ,\; k\ge 1$,

\vskip2.5cm 
\hskip2.5cm
\setlength{\unitlength}{0.1cm}
\begin{picture}(120,60)
  \epsfxsize=12.0cm
  \epsffile{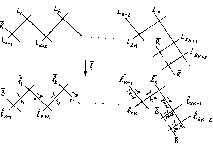}
\end{picture}

\begin{center}
$Z=3L_1 +5L_2+\ldots +(2k+1)L_k +
2L_{k+1} +4L_{k+2}+\ldots +2kL_{2k} +(k+1)L_{2k+1} + (k+1)L_{2k+2} \, ;$
\end{center}

\newpage
\noindent 4) The singularity  $D_{2k+3}: x(y^2+x^{2k+1}) ,\; k\ge 1$,

\vskip3cm 
\hskip2.5cm
\setlength{\unitlength}{0.1cm}
\begin{picture}(120,60)
  \epsfxsize=12.0cm
  \epsffile{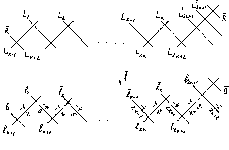}
\end{picture}

\begin{center}
$Z=3L_1 +5L_2+\ldots +(2k+1)L_k +
2L_{k+1} +4L_{k+2}+\ldots +2kL_{2k}  + (2k+2)L_{2k+2}
+(k+1)L_{2k+1} \, ; $
\end{center}
\noindent £¤¥ $L_{2k+1} =L_{2k+1}' +L_{2k+1}'' \, ;$

\vspace{0.1cm}

\noindent 5) The singularity  $E_6: x^3+y^4 \, ,$

\vskip3cm 
\hskip2cm
\setlength{\unitlength}{0.1cm}
\begin{picture}(80,40)
  \epsfxsize=7.0cm
  \epsffile{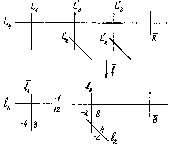}
\end{picture}

\begin{center}
$Z=3L_1 +2L_2+4L_3 +6L_4 \, ;$
\end{center}

\newpage \noindent 6) The singularity  $E_7: x(x^2+y^3) \, ,$

\vskip3cm 
\hskip2cm
\setlength{\unitlength}{0.1cm}
\begin{picture}(80,40)
  \epsfxsize=8.0cm
  \epsffile{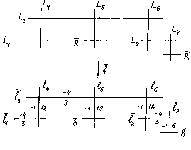}
\end{picture}

\begin{center}
$Z=3L_1 +5L_2+9L_3 +6L_4 +5L_5 +8L_6 +3L_7  \, ;$
\end{center}

\noindent 7) The singularity  $E_8: x^3+y^5 \, ,$

\vskip3cm 
\hskip2cm
\setlength{\unitlength}{0.1cm}
\begin{picture}(120,60)
  \epsfxsize=10.0cm
  \epsffile{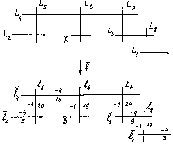}
\end{picture}

\begin{center}
$Z=3L_1 +5L_2+9L_3 +15L_4 +10L_5 +8L_6 +12L_7 +6L_8  \, .
\hphantom{aaa} \rule{6pt}{6pt}$
\end{center}
\end{lem}

\vspace{0.1cm}

{\bf 3.4.} {\it Defect of a singularity.} \label{s:3.4}
Define a {\it defect} $\delta$ of a A-D-E-singularity
by the formula
$$
\delta = \frac12 (\bar R \cdot Z) \, .
$$

\begin{cor}
For different types of A-D-E-singularities the defect equals
$$
\delta =
\left\{
\begin{array}{ll}
\left[ \frac{n+1}{2} \right] & \mbox{for type\hphantom{a}} A_n ; \\
\left[ \frac{n}{2} \right] +1 & \mbox{for type\hphantom{a}} D_n ; \\
\left[ \frac{n+1}{2} \right] & \mbox{for types\hphantom{a}} E_n ,\; n=6,7,8 \, .  \\
\end{array}
\right.
$$

\noindent In particular, for the type $A_1$ (nodes) and $A_2$ (cusps) the defect
$\delta =1$.

\vspace{0.2cm}
Actually one can show that defect $\delta$ is the
$\delta$-invariant (genus) of the one-dimensional A-D-E-singularity.

\end{cor}

\section{Numerical invariants of a generic covering}

   Now we consider a global situation. Let $X$ be a surface with
A-D-E-singularities,
$$
Sing\;X = \sum_{k\ge 1}a_kA_k +\sum_{k\ge 4}d_kD_k + \sum_{k=6,7,8}e_kE_k \, ,
$$
that means that $X$ has $a_k$ singularities of type $A_k$, $d_k$ -- of type
$D_k$ and $e_k$ -- of type $E_k$. Let $f: X \to {\bf P}^2$ be a generic
covering of degree $N$, and $B\subset {\PP}^2$ be the discriminant curve.
Let $\deg B =d$ and let $B$ has $n$ nodes and $c$ cusps,
$n_s =a_1$ and $c_s =a_2$ of which originates from $Sing\:X$, and $n_p$
and $c_p$ are p-nodes and p-cusps.
Let $R\subset X$ be the ramification curve, $f^*(B)= 2R+C$,
and $L\subset X$ be the preimage of a generic line $l\subset {\PP}^2$.
Let $\pi : S \to X $ be the minimal resolution of $X$, and
$\bar f = f \circ \pi : S \to {\PP}^2$.
Denote by $\bar R$ and $\bar L$ the proper transforms  of $R$ and $L$ on $S$.
Then $\bar R$ is a normalization of the curve $R\simeq B$, and
$\bar L \simeq L$.

\vspace{0.2cm}

{\bf 4.1.} {\it The canonical class $K_S$ and the canonical cycle $Z$.}
Let
$$
Z=\sum_{x\in Sing X}Z_x
$$
be the canonical cycle of $S$, where $Z_x$ are the canonical cycles of
singularities $x \in Sing\: X$.
It follows from 3.2 that
$$
K_S = (f\circ \pi )^*K_{{\PP}^2} +\bar R +Z =
-3\bar L + \bar R + Z\, .
\eqno (4.1)
$$
Besides, the singularities of $X$ being Gorenstein, the divisor $R$
is locally principal, and
$$
{\pi}^*(R) =\bar R +Z \, .
\eqno (4.2)
$$

\vspace{0.2cm}

{\bf 4.2.} {\it The intersection numbers.}

\begin{lem}
The intersection numbers of $\bar L ,$ $\bar R$ and $Z$ on S
are equal
$$
\left( \bar L^2 \right) = N \; , \; \bar L \cdot \bar R = d \; , \;
\bar L \cdot Z=0 \; , \bar R \cdot Z =2{\delta}_X \, , \;
\left( Z^2 \right) =-2{\delta}_X \: ,
\eqno (4.3)
$$
where
$$
{\delta}_X = \sum_{x\in Sing X}{\delta}_x =
\sum a_k\left[ \frac{k+1}{2} \right] +
\sum d_k \left(\left[ \frac{k}{2} \right] +1 \right) +
\sum e_k\left[ \frac{k+1}{2} \right]
\eqno (4.4)
$$
is the defect of the surface $X.$
\end{lem}

\dvo
Obviously, we have
$\left( \bar L^2 \right) =deg\,f = N , $  and $ \; \bar L\cdot\bar R =
deg\,B = d$.
By 3.4 we have $\bar R \cdot Z =2{\delta}_X$.
The divisor $Z$ being exceptional, we have
$\pi (Z) =Sing\: X ,\: \dim \pi (Z) =0 $, and
$\bar L ={\pi}^*(L)$, $\bar R +Z ={\pi}^*(R)$, and therefore,
$\bar L \cdot Z=0$, and $(\bar R +Z) \cdot Z=0 $, and, consequently,
$\left( Z^2 \right) =-(\bar R \cdot Z) $. $\rule{6pt}{6pt}$

It remains to compute $\left( \bar R^2 \right) $.

\vspace{0.1cm}

{\bf 4.3.} {\it The evenness of degree $deg\,B =d =2\bar d$.}
The restriction of $\bar f$ to $\bar L $ ,
$\bar L \to l\simeq {\PP}^1 $, is a covering of degree $N$, with
ramification indices 2 at the points of intersection of $\bar L$ and
$\bar R$. We have $\bar L \cdot \bar R =d $, and from Hurwitz formula
we obtain $2g(\bar L) -2 =-2N+d$. It follows that $\deg B =d$
is even. Let $d=2\bar d $. Besides, since
$$
g(\bar L) =\frac12 d +1-N \ge 0 \, ,
$$
we obtain a bound for the degree of covering,
$$
N\le \bar d +1 \, .
$$

\vspace{0.1cm}

{\bf 4.4.} {\it The self-intersection number $\left( \bar R^2 \right) $ and the
arithmetical genus of the curve $R$.}
Denote by $\delta$ the defect of the curve $B$,
$$
{\delta} ={\delta}_B = \sum_{s\in Sing B}{\delta}_s =
n+c+{\delta}_0 ,
\eqno (4.5)
$$
\noindent where
$$
{\delta}_0 = \sum_{x\in Sing B , \; x \;
\mbox{\footnotesize not} \; A_1 \mbox{\footnotesize and} A_2}{\delta}_x .
\eqno (4.6)
$$
\noindent The numbers $\delta$ and ${\delta}_0$ are the extremal values of
defects ${\delta}_X$ of surfaces $X$ with given discriminant curve $B:$
${\delta}_0$ corresponds to a surface $X$, all nodes and cusps of which are
p-nodes and p-cusps, $n=n_p ,$ $c=c_p$, and $\delta$
corresponds to a surface $X $ (for example, a 2-sheeted covering of
$\PP^2$), all nodes and cusps of which are s-nodes and s-cusps,
$n=n_s ,$ $c=c_s$.

At first we express the geometric genus of $B ,$ $g=g(B)=g(\bar R)$,
in terms of the defect $\delta$. For this we consider a surface $X$,
which is a 2-sheeted covering of ${\PP}^2$ with the discriminant curve $B$.
In this case $\left( Z^2 \right) =-(\bar R \cdot Z) =-2{\delta} $, and
$f^*(B) =2R $ and, consequently, $d\cdot \bar L \sim 2\bar R +2Z $, because
$B\sim d\cdot l $. From (4.1) and the adjunction formula
$g(\bar R) =\frac{\displaystyle (\bar R,\bar R +K_{\bar X})}{\displaystyle2}
+1 $
we obtain
$$
g=\frac{(d-1)(d-2)}{2} -{\delta} \, .
\eqno (4.7)
$$

\vspace{0.2cm}
If it is known that the defect $\delta $ coincides with the
$\delta$-invariant of a one-dimensional singularity, then this formula
coincides with the well known formula for the geometric genus
$g(R) \stackrel{df}{=}g(\bar R)$ of a singular curve
$R$, $g(R) =p_a(R) -\sum_{x\in Sing R}{\delta}_x $.

We return to a generic covering $X$ of degree $N ,$ $n=n_s +n_p $,
$c=c_s +c_p $. Then
$$
{\delta}_X = n_s +c_s + {\delta}_0 =\delta -n_p -c_p .
\eqno (4.8)
$$

\begin{lem}
The self-intersection number of the proper transform
of the ramification curve $\bar R \subset S$ is equal
$$
\left( \bar R^2 \right) =3\bar d +g-1- {\delta}_X \, ,
\eqno (4.9)
$$
and
$$
\left( \bar R +Z \right)^2 =3\bar d +g-1+ {\delta}_X
=3\bar d +p_a(R) -1 \, ,
\eqno (4.10)
$$
where
$$
p_a(R) =g +{\delta}_X =\frac{\displaystyle (d-1)(d-2)}{\displaystyle 2}
-n_p -c_p
\eqno (4.11)
$$
is the arithmetical genus of $R$.
\end{lem}

\dvo
From (4.1) and the adjunction formula
$2g(\bar R) -2 = (\bar R,\bar R +K_S) = (\bar R,-3\bar L +2\bar R +Z)$
we obtain
$\left( \bar R^2 \right) =\frac32 (\bar R \cdot \bar L) +g -1 -
\frac12 (\bar R \cdot Z) $. Applying formulae (4.3), we obtain the proof.
$\rule{6pt}{6pt}$

From formulae (4.1), (4.3) and (4.9) we obtain a corollary.

\begin{cor}
$$
\left( K_S^2 \right) =9N -9\bar d +p_a(R) -1 ,
\eqno (4.12)
$$
or, substituting $p_a(R)$ from (4.11),
$$
\left( K_S^2 \right) =9N  +\frac12 d(d-12) -n_p -c_p .
\eqno (4.12')
$$
\end{cor}

\vspace{0.1cm}

{\bf 4.5.} {\it A bound for the covering degree.}
\begin{lem}
For a generic covering of degree $N$ with discriminant
curve of degree $d=2\bar d $ and genus $g$, we have
$$
N \le \frac{4\bar d^2}{3\bar d +g-1+{\delta}_X} \, ,
\eqno (4.13)
$$
where ${\delta}_X$ is the defect of singularities of $X$, and moreover,
the equality holds if and only if $\bar L \equiv mK_S $
for some $m\in {\QQ}^* $, or $mK_S \equiv 0$.
\end{lem}

\dvo
Applying Hodge index theorem to divisors
$\bar L$ and ${\pi}^*(R) =\bar R +Z $ on $S$, we obtain
$$
\left|
\begin{array}{cc}
\bar L^2 & (\bar L,\bar R +Z)\\
(\bar L,\bar R +Z)& (\bar R +Z)^2
\end{array}
\right|
=
\left|
\begin{array}{cc}
N & d\\
d & 3\bar d +g-1+ {\delta}_X
\end{array}
\right|
\le 0 \, ,
$$
and it is the desired inequality. The equality holds only if $\bar L$
and $\bar R +Z $ are linear dependent in the N\'{e}ron-Severi group
$NS(\bar X)\otimes{\bf Q} $. Since $K_S = -3\bar L + \bar R + Z $,
we obtain the assertion about possible equality. $\rule{6pt}{6pt}$

\vspace{0.1cm}

{\bf 4.6.} {\it The topological Euler characteristic $e(S).$}

\begin{lem}
The topological Euler characteristic of a surface $S$,
obtained by the minimal resolution of singularities of $X$, is connected
with the defect ${\delta}_X$ and invariants of a generic covering $f$ by
a formula
$$
e(S) = 3N+2g-2+2{\delta}_X -c_p ,
\eqno (4.14)
$$
where $N=\deg f $, and $c_p$ is the number of p-cusps on $B$
(or the number of pleats of $f$).
\end{lem}

{\it Proof} is obtained in the same way as in the case of a non-singular
surface $X$ ([K], $\S 1$ Lemma 7), considering a generic pencil of lines on
${\PP}^2$ and the corresponding hyperplane sections on $S$, and lifting the
morphism $\bar f: S \to {\PP}^2$ to a morphism of fiberings of curves
over ${\PP}^1.$ One can obtain a proof by direct computations. At first
we find $e(X) = 3N-e(B)-n_p -c_p $ by
considering the finite covering $f: X \to {\PP}^2 $, the stratification
${\PP}^2 =({\PP}^2\setminus B)\cup (B\setminus Sing\; B)\cup Sing\; B$,
and applying the additivity property of Euler characteristic, and then we
find $e(S)$. $\rule{6pt}{6pt}$

From Noether's formula $(K_S^2) +e(S) =12p_a$ and formulae (4.12) and (4.14)
we have $12p_a =12N -9\bar d +3p_a(R) -3 -c_p $.
Substituting $p_a(R)$ from (4.11), we obtain a corollary.

\begin{cor} The Euler characteristic of the structure sheaf ${\cal O}_S$
equals
$$
p_a =1-q-p_g =N +\frac{\bar d (\bar d -3)}{2} -\frac{n_p}{4} -\frac{c_p}{3} .
\eqno (4.15)
$$
\end{cor}

Thus, as in the case of a non-singular surface $X ,$ we obtain

\begin{cor}
$$
n_p \equiv 0\;(mod\;4) \: , \: c_p \equiv 0\;(mod\;3) .
$$
\end{cor}

\section{Proof of the main inequality.}

{\bf 5.1.} {\it A fiber product of two generic coverings.}
Let a curve $B$ be a common discriminant curve for two generic coverings
$f_1:X_1\to {\PP}^2$ and $f_2:X_2\to {\PP}^2$ of degrees
$\deg f_1=N_1$ and $deg f_2= N_2$. Let
$$
Sing\; B = nA_1 +cA_2 +
\sum_{k>2}a_kA_k +\sum_{k\ge 4}d_kD_k + \sum_{k=6,7,8}e_kE_k \; .
$$
With respect to a pair of coverings $f_1$ and $f_2$ nodes and cusps of $B$
are subdivided into four types,
$$ n=n_{ss}+n_{sp}+n_{ps}+n_{pp} \; , \quad
c=c_{ss}+c_{sp}+c_{ps}+c_{pp} ,
\eqno (5.1)
$$
where $n_{\flat \sharp}$ and $c_{\flat \sharp}$ are numbers of
$\flat \sharp$-nodes and $\flat \sharp$-cusps of $'$. In particular,
$n_{ss} +n_{sp} =a_1 $ is the number of singularities of type $A_1$,
and $á_{ss} +á_{sp} =a_2 $ is the number of singularities of type
$A_2$ on the surface $X_1$.

Consider a normalization $X$ of the fiber product
$X^{\scriptstyle \times} =X_1{\times}_{{\bf P}^2}X_2$
and the corresponding commutative diagram
\\

\unitlength=1.00mm
\special{em:linewidth 0.4pt}
\linethickness{0.4pt}
\begin{picture}(72.00,52.00)
\put(60.00,10.00){\makebox(0,0)[ct]{$(x,y)\in \mathbb C^2$}}
\put(80.00,50.00){\makebox(0,0)[cb]
{$X=\mathbb C^2 \ni (z_1,z_2) $}}
\put(55.00,30.00){\makebox(0,0)[rc]{$ (z_1,y)\in \mathbb C^2=X_1$}}
\put(78.00,30.00){\makebox(0,0)[lc]{$X_2= \mathbb C^2 \ni (x,z_2)\:$ ,}}
\put(66.00,48.00){\vector(0,-1){15.00}}
\put(71.00,40.00){\makebox(0,0)[rc]{$g_0 $}}
\put(66.00,28.00){\vector(0,-1){15.00}}
\put(71.00,22.00){\makebox(0,0)[rc]{$f_0 $}}
\put(64.00,48.00){\vector(-3,-4){10.67}}
\put(68.00,48.00){\vector(3,-4){10.67}}
\put(53.00,27.00){\vector(3,-4){10.67}}
\put(79.00,27.00){\vector(-3,-4){10.67}}
\put(57.00,20.00){\makebox(0,0)[rc]{$f_1 $}}
\put(75.00,20.00){\makebox(0,0)[lc]{$f_2 $}}
\put(58.00,42.00){\makebox(0,0)[rc]{$g_1$}}
\put(64.00,30.00){\makebox(0,0)[lc]{$X_0$}}
\put(74.50,42.00){\makebox(0,0)[lc]{$g_2$}}
\put(160.00,30.00){\makebox(0,0)[rc]{$(*_2)$}}
\end{picture}

\noindent The surface $X$ is a $N_1N_2$-sheeted covering of ${\PP}^2$
and it has at most A-D-E-singularities, which lie over $Sing\;B $.

\noindent {\bf Lemma.} If coverings $f_1$ and $f_2$ are non equivalent, then the
surface $X$ is irreducible.

{\it Proof} is word for word the same as in the case of generic coverings
of non-singular surfaces (\cite{K} Proposition 2). $\rule{6pt}{6pt}$

We set
$$
g_1^{-1}(R_1) =R+C  ,
$$
where $R$ is a part, which is mapped by $g_2$ onto $R_2 $,
and $C$ is a part, which is mapped $g_2$ onto $C_2 $.
We are interested in the intersection number of $R$ and $C$ after
a resolution of singularities of $X$ in a neighbourhood of the curve $R+C $.

Consider a restriction $R+C \to R_1$ of the covering $g_1$ over the curve
$R_1 $. As follows from 2.2.1 and 2.2.2, it is an \'{e}tale covering
of degree $N_2$ over a generic point $x_1 \in R_1$, where  $R \to R_1$ is
a 2-sheeted, and $C \to R_1$ is a $(N_2-2)$-sheeted covering. The same
picture is over a point $x_1 \in R_1 $, which is a s-point of $X_1 $,
lying over a ss-point of $B $.

Denote by $\tilde{\pi}: S \to X $ a minimal resolution of singularities of
$X $, and denote by $\tilde R$ and $\tilde '$ the proper transforms of $R$
and $C$ on $S $.
Our goal is to calculate the intersection numbers $({\tilde R}^2) $,
$(\tilde R \cdot \tilde C )$ and $({\tilde C}^2) $, and also the analogous
intersection numbers for divisors $\tilde{\pi}^{-1}(R) =\tilde R +Z_R$ and
$\tilde{\pi}^{-1}(') =\tilde ' +Z_' $, where $Z_R$ and $Z_C$ are the
sums of canonical cycles corresponding to singular points
$x \in Sing\; X $ and lying on $R$ and $C$ respectively.

\vspace{0.1cm}

{\bf 5.2.} {\it The structure of a fibre product over a neighbourhood of a
singular point of the discriminant curve.}
Let $U\subset {\PP}^2$ be a sufficiently small neighbourhood (in complex
topology) of a point $b \in Sing\; B $. The preimage $f_1^{-1}(U)$ is a
disjoint union of two parts, $f_1^{-1}(U) = V_1 \sqcup V_1' ,$ where
$V_1$ is a part containing the ramification curve $R_1 $, and $V_1'$ is
a part not containing $R_1$ and \'{e}tale mapped to $U $.
Analogously $f_2^{-1}(U) = V_2 \sqcup V_2'$. Then $f^{-1}(U)$
is a disjoint union of four open sets -- of normalizations of fibre products
$W= \overline{V_1{\times}_UV_2} ,$  $W'=  \overline{V_1{\times}_UV_2'} ,$
$\overline{V_1'{\times}_UV_2}$ and $\overline{V_1'{\times}_UV_2'} $.
And only $W$ and $W'$ meet the curve $g_1^{-1}(R_1) $.
The open sets $W\subset X$ were studied in detail in $\S 2 $. The surface
$X$ in the neighbourhood $W$ is non-singular except the case of
ss-points $b$. The open set $W'$ consists of $N_2-k $ components
($k=2,3$ or $4$ depending on the type of the singular point $b$),
which are mapped isomorphically onto $V_1$. And $W'$ does not meet
$R $, and $W' \cap C$ consists of $N_2 -k$ components isomorphic to
$V_1 \cap R_1 $.

It follows from the investigation of the local structure of $X$ in $\S 2$
that $X$ and the curves $R$ and $C$ are of the following form over
neighbourhoods of singular points $b \in Sing\; B $ of different types.

1) Over a ss-point $b$ the neighbourhood $W$ has 2 , and $W'$ has $N_2-2$
components, which are mapped isomorphically onto $V_1$ by the map $g_1 $.
Correspondingly $R\cap W$ consists of two, and $C\cap W'$ consists of
$(N_2-2)$ components isomorphic to $R_1 \cap V_1 $.

\vspace{0.1cm}

2) Over a sp-point $b\in B $ of type $A_1$ the neighbourhood $W'$ consists of
$(N_2-4)$ components isomorphic to $V_1$ and having a singular point of type
$A_1 $. Correspondingly $C$ consists of $N_2-4$ nodal curves.
The neighbourhood $W$ consists of two components: see Fig. 4, where
$$
R = R^{\prime ,\prime} + R^{\prime \prime ,\prime} \; , \: \mbox{and \quad} \:
C = R^{\prime ,\prime \prime} + R^{\prime \prime ,\prime \prime}
$$
\noindent (it ought to change places of the left and right parts of Fig. 4,
$g_1$ stands for $g_2 $, and $g_2$ -- for $g_1 $ ). We see that in
the neighbourhood $W$ the curves $R$ and $C$ are non-singular and meet
transversally in two points.

\vspace{0.1cm}

3) Over a ps-point $b\in B $ of type $A_1$ the neighbourhood
$V_1 \subset X_1$
consists of two components and on each of them the map $f_1$ has
a fold. The neighbourhood $W'$ consists of disjoint union of $(N_2-2)$
pieces isomorphic to $V_1 $. The neighbourhood $W$ consists of two components:
see Fig. 4, on which
$$
R = R^{\prime ,\prime} + R^{\prime \prime ,\prime} \; , \: \mbox{and \quad} \:
C = \emptyset .
$$
We see that on $W$ the curve $R$ is non-singular and does not meet
$C $.

\vspace{0.1cm}

4) Over a pp-point $b\in B $ of type $A_1$ the neighbourhood
$V_1 \subset X_1$
consists of two irreducible components and on each of them $f_1$ has a fold.
The neighbourhood $W'$ is non-singular and consists of
$N_2-4 $ components isomorphic to $V_1 $. The neighbourhood $W$
is represented on Fig. 5, on which
$$
R = R^{\prime ,\prime} + R^{\prime ,\prime \prime} +
R^{\prime \prime ,\prime} + R^{\prime \prime ,\prime \prime} , \;
C = C^{\prime ,\prime} + C^{\prime \prime ,\prime}  .
$$
We see that the curves $R$ and $C$ are non-singular and do not meet.

\vspace{0.1cm}

5) Over a sp-point $b\in B $ of type $A_2$ the neighbourhood $V_1$ has
a singular point of type $A_2 $, and $W'$ consists of $(N_2-3)$ components
isomorphic to $V_1 $. The neighbourhood $W$ is pictured on Fig. 7, on which
$$
R = L_2 +L_3 \; , \: C = L_1  .
$$
We see that $R$ has a double point, $C$ is non-singular and intersect
transversally each of the branches of $R$ at the intersection point, and,
consequently, $(R\cdot C) =2 $.

\vspace{0.1cm}

6) Over a ps-point $b\in B $ of type $A_2$ the neighbourhood $V_1$ is
non-singular, and $W'$ consists of $(N_2-2)$ components isomorphic to $V_1 $.
The neighbourhood $W$ is pictured on Fig. 7 (on which it ought to change
places of the left and right parts, $g_1$ stands for $g_2 $, and
$g_2$ -- for $g_3 $ ), where
$$
R = L_2 +L_3 \; , \:
C = \emptyset \; .
$$
We  see that $R$ has a double point and does not meet $C$.

\vspace{0.1cm}

7) Over a pp-point $b\in B $ of type $A_2$ the neighbourhood $W'$ consists of
$N_2-3 $ components isomorphic to $V_1 $. The neighbourhood $W$
is pictured on Fig. 8, on which
$$
R = R_3 +L_3 \; , \; C=L_2 \; .
$$
We see that $R$ is non-singular and meets with $C$ transversally
at one point.

From the obtained local description it follows that the surface $X$ is
non-singular at the points of intersection of $R$ and $C$,
and the intersection is not void only over the points $b\in B$
of types: over sp-points of type $A_1 $, where $(R\cdot C)=2 $,
over sp-points of type $A_2 $, where $(R\cdot C)=2 $,
and over pp-pointe of type $A_2 $, where $(R\cdot C)=1 $. Therefore,
$$
(\tilde R\cdot \tilde C) =2n_{sp} +2c_{sp} +c_{pp}.
\eqno (5.3)
$$

\vspace{0.1cm}

{\bf 5.3.} {\it A lift of the fibre product to a resolution of the
discriminant curve.}
To compute intersection numbers on $S$ we consider firstly an auxiliary
surface $\bar X $, which is not a minimal resolution of $X $, and then we
\lz descend' to $S $. Let $\sigma : \bar {\PP}^2 \to {\PP}^2 $
be a composition of $\sigma$-processes resolving the curve $B$ and needed to
obtain a minimal resolution of a double plane singularities, lying over $B$
(see $\S 3$), and, besides, let $\sigma$ includes two additional
$\sigma$-processes as in 2.4.5 for each cusp, which is not a ss-cusp.
Consider a lift of the diagram $(*_1)$ to $\bar {\PP}^2 $, namely consider the
diagram

$$\xymatrix
{
& & S \ar[ld]_{\widetilde \pi }  & & \\
 & X \ar[ld]_{g_1} \ar@{-->}[rdd]^(.30){g_2} & & \overline{X} \ar[ll]_{\pi }
 \ar[lu]_{\overline \pi }
\ar[ld]_{\overline g_1} \ar[rdd]^{\overline g_2} &   \\
X_1 \ar[rdd]_{f_1} & & \overline X_1 \ar[ll]^{\pi_1} \ar[rdd]^(.30){\overline f_1} & & \\
& & X_2 \ar@{-->}[ld]_{f_2} & & \overline X_2 \ar@{-->}[ll]_{\pi _2}
\ar[ld]^{\overline f_2} \\
& \mathbb P^2 & & \overline{\mathbb P}^2 \ar[ll]_{\sigma } &
}
\eqno (5.4)
$$

\noindent in which $\bar X_i$ and $\bar X$ are normalizations of
$X_i\times_{{\PP}^2}{\bar {\PP}}^2$ and
$X\times_{{\PP}^2}{\bar {\PP}}^2 $.
Then morphisms \lz on the right wall' of diagram (5.4) are finite coverings.
The surface $\bar X$ is non-singular, and
$\bar \pi : \bar X \to S $ blows down the \lz superfluous' exceptional curves
of the first kind. Let $\bar R_1$ be the proper transform of $R_1$ on
$\bar X_1 $, and $\bar R$ and $\bar C$ (respectively $\tilde R$ and
$\tilde C$) be the proper transforms of $R$ and $C$ on $\bar X$
(respectively on $S$).
Then $\bar g^*_1(\bar R_1) =\bar R +\bar C $, and $\bar R \to \bar R_1$ and
$\bar C \to \bar R_1$ are finite coverings of degree 2 and $N_2 -2$
respectively, and $\bar R$ and $\bar C$ are disjoint. Therefore,
$$
\left( \bar R^2\right) =2\left( \bar R_1^2\right)  , \:
\left( \bar C^2\right) =(N_2-2)\left( \bar R_1^2\right)  , \:
\bar R \cdot \bar C =0  .
\eqno (5.5)
$$

Actually from 3) and 4) one can see that over ps- and pp-nodes $b$
in a neighbourhood of $R+C$ the surface $X$ is non-singular, the curves
$R$ and $C$ are non-singular and disjoint. Therefore, one can suppose that
ps- and pp-nodes on $B$ are not blown up (and on the surface $S$ there
remain singular points, which lie over these nodes).

\vspace{0.1cm}

{\bf 5.4.} {\it Computation of intersection numbers.}
First we find $({\bar R_1}^2) $. Recall that by (4.9) we have on the
minimal resolution $\tilde X_1$ of the surface $X_1$
$$
({\tilde R_1}^2) =3\bar d + g-1-{\delta}_1 ,
\eqno (5.6)
$$
\noindent where ${\delta}_1 ={\delta}_{X_1} =n_s +c_s +{\delta}_0 $, and
$n_s =n_{ss} +n_{sp}$ and $á_s =á_{ss} +á_{sp}$ are the numbers of singular
points of type $A_1$ and $A_2$ on the surface $X_1 $.

Let ${\pi}_1=\tilde {\pi}_1 \circ \bar {\pi}_1 $, where
$\tilde {\pi}_1: \tilde X_1 \to X_1$ is a minimal resolution, and
$\bar {\pi}_1: \bar X_1 \to \tilde X_1$ is the blowing down of the
"superfluous" exceptional curves. The surface $\bar X_1$ differs from the
surface $\tilde X_1$ only over the cusps of $B $, which are not ss-cusps.
Let $\bar U ={\sigma}^{-1}(U) $, and $\bar V_1 ={\pi}_1^{-1}(V_1) ,
\tilde V_1 ={\tilde {\pi}_1}^{-1}(V_1) $ be neighbourhoods of $\bar X_1$
and $\tilde X_1$ lying over $\bar U$ and containing the proper transform
of $R_1 $. Analogously $\bar V_1'$ and $\tilde V_1' $.

For a sp-cusp $b\in B$ the blowing down
$\bar V_1 \stackrel{\bar{\pi}_1}{\longrightarrow}{\tilde V_1}$
is represented on Fig. 12. We see that $\bar{\pi}$ includes one $\sigma$-process
with a centre $R_1 $. For ps- and pp-cusps $b\in B$ the blowing down
$\bar V_1 \stackrel{\bar{\pi}_1}{\longrightarrow}{\tilde V_1}$
is represented on Fig. 13 (where $R$ stands for $R_1 $, and $C$ stands for
$C_1 $). We see that it needs two $\sigma$-process with centres on
$R_1$ to disjoint $C_1$ and $R_1 $. Therefore,
$$
({\bar R_1}^2) =({\tilde R_1}^2) -c_{sp} -2c_{ps} -2c_{pp} .
\eqno (5.7)
$$

\vspace{0.1cm}

   Now we examine how the intersection numbers
$({\bar R}^2)$ and $({\bar C}^2)$ change under the blowing down $\bar \pi $.
For a neighbourhood $U\subset {\PP}^2$ of a point $b\in Sing\: B$ set
$\bar W ={\pi}^{-1}(W) , \; \bar W' ={\pi}^{-1}(W') , \;
\tilde W ={\tilde \pi}^{-1}(W) , \;
\tilde W' ={\tilde \pi}^{-1}(W') $. Then
$\bar g_1^{-1}(\bar V_1) =\bar W \sqcup \bar W' $. We examine one after
another the blowing down $\bar{\pi}: \bar X \to S $ in neighbourhoods
$\bar W \sqcup \bar W' \subset \bar X$ separately for different types of
singular points $b\in Sing\; B$  (the numbering of cases corresponds to the
numbering of cases in 5.2).

2) For a sp-point $b$ of type $A_1$ the neighbourhood $\bar W' $
is a disjoint union of $(N_2-4)$ open sets isomorphic to $\bar V_1 $ --
to the minimal resolution of singular points of type $A_1 $. The neighbourhood
$W$ is represented on Fig. 4 , and $\bar {\pi}: \bar W \to W $ is a blowing up
of two points $R^{\prime ,\prime} \cap R^{\prime ,\prime \prime} $ and
$R^{\prime \prime ,\prime} \cap R^{\prime \prime ,\prime \prime} $.
Therefore, the blowing down $\bar {\pi}: \bar W \to \tilde W \simeq W $
increases $(\bar R^2)$ and $(\bar C^2)$ on 2 for one point $b$ and,
consequently, on $2n_{sp}$ for all points of this type.

5) For a sp-point $b$ of type $A_2$ the neighbourhood $\bar W $
is represented on the upper part of Fig. 10 . It is obtained from the
neighbourhood $W $, pictured on Fig. 7 , by blowing up the point of
intersection of lines $L_1 , L_2$ and $L_3 $, and then by blowing up 5 points
on the glued line $E_{6,3}$ and not lying on the proper transform of these
lines. The blowing down $\bar {\pi}: \bar W \to \tilde W \simeq W $ is
the converse procedure, i.e. the blowing down of five exceptional curves of
the first kind, and then blowing down the curve $E_{6,3} $. In this case
$R=L_2 +L_3 $, and $C=L_1 $. Since
$(R^2) =(L_2^2) +2(L_2\cdot L_3) +(L_3^2) $ and $(L_2^2) ,$
$(L_3^2)$ are diminished on 1, and $L_2$ and $L_3$ are no longer intersected
after the $\sigma$-process with the centre at the point $L_2 \cap L_3 $,
the blowing down $\bar {\pi}$ increases $(\bar R^2)$ on 4 for one point $b$
and on $4c_{sp}$ for all points of this type.

The neighbourhood $\bar W'$ consists of $(N_2-3)$ components isomorphic to
$\bar V_1 $, for each of which $\bar {\pi}$ is represented on Fig. 12 .
As above for $(\bar R_1^2)$, we see that the blowing down $\bar {\pi}$ increases
$(\bar C^2)$ on $(N_2 -3)+1$ (taking account of the neibourhood $\bar W$ )
for one point $b$ and on $(N_2 -2)c_{sp}$ for all points of this type.

6) For a ps-point $b$ of type $A_2$ the neighbourhood $\bar W $ and the
blowing down $\bar {\pi}: \bar W \to \tilde W \simeq W $ are the same
as in 5) , but in this case $R=L_2 +L_3 $, and $C \cap W =\emptyset $.
Therefore, as in 5) we obtain that the blowing down $\bar {\pi}$ increases
$(\bar R^2)$ on $4c_{ps} $.

The neighbourhood $\bar W'$ consists of $(N_2-2)$ components isomorphic to
$\bar V_1 $, for each of which $\bar {\pi}$ is represented on Fig. 13 .
As above for $(\bar R_1^2)$, we see that the blowing down $\bar {\pi}$ increases
$(\bar C^2)$ on $2(N_2 -2)$ for one point $b$ and on $2(N_2 -2)c_{ps}$
for all points of this type.

7) For a pp-point $b$ of type $A_2$ the neighbourhood $\bar W $
consists of two components: one is the same as in 5) and the other is the
same as $\bar V_1'$ and represented on the left side of Fig. 13 .
Since in the neibourhood $W $,  represented on Fig. 8,
$R=R_3 +L_3 $, and $C=L_2 $, we obtain that the blowing down
$\bar {\pi}: \bar W \to \tilde W \simeq W $ increases $(\bar R^2)$
on $1+2=3$ for one point  $b$ and on $3c_{pp}$ for all points of this type.
Besides, $(\bar C^2)$ is increased on $c_{pp} $.

The neighbourhood $\bar W'$ consists of $(N_2-3)$ components isomorphic to
$\bar V_1' $, and is represented on Fig. 13 (on which $C$ stands for $R$ ).
Therefore, taking account of the neighbourhood $\bar W $,
the blowing down $\bar {\pi}$ increases $(\bar C^2)$ on
$2(N_2 -3)c_{pp} + c_{pp} =(2N_2-5)c_{pp} $.

Summing all modifications of $(\bar R^2)$ and $(\bar C^2) $, we obtain
$$
(\tilde R^2) =(\bar R^2) + 2n_{sp} +4c_{sp} +4c_{ps} + 3c_{pp} ,
\eqno (5.8)
$$
$$
(\tilde C^2) =(\bar C^2) + 2n_{sp} +(N_2-2)c_{sp} +
2(N_2-2)c_{ps} + (2N_2-5)c_{pp} .
\eqno (5.9)
$$
Applying (5.5) and substituting $({\bar R_1}^2)$ from (5.7), we obtain
$$
(\tilde R^2) =2(({\tilde R_1}^2)-c_{sp} -2c_{ps} -2c_{pp})
+ 2n_{sp} +4c_{sp} +4c_{ps} + 3c_{pp} =
$$
$$
= 2({\tilde R_1}^2) + 2n_{sp} +2c_{sp} - c_{pp} ,
\eqno (5.10)
$$
$$
(\tilde C^2) =(N_2-2)(({\tilde R_1}^2)-c_{sp} -2c_{ps} -2c_{pp}) +
2n_{sp} +(N_2-2)c_{sp} + 2(N_2-2)c_{ps} + (2N_2-5)c_{pp} =
$$
$$
= (N_2-2)({\tilde R_1}^2) +
2n_{sp} - c_{pp} .
\eqno (5.11)
$$

\vspace{0.1cm}

{\bf 5.5.} {\it Computation of intersection numbers (continuation).}
Now we find $(\tilde R +Z_R)^2 ,$ $(\tilde C +Z_C)^2 $ and
$(\tilde R +Z_R)\cdot(\tilde C +Z_C) $, where the divisor $Z_R $,
respectively $Z_C $, equals to $\sum Z_x $, where $Z_x$ is the canonical
cycle of a point $x\in Sing\; X $, and the summation runs over $x\in R $,
respectively $x\in C .$
The analogous sums $\sum {\delta}_x$ we denote by ${\delta}_R$ and
${\delta}_C$ respectively. By (4.2) we have
$$
(\tilde R\cdot Z_R) =-(Z_R^2) =2{\delta}_R \; ,
(\tilde C\cdot Z_C) =-(Z_C^2) =2{\delta}_C \; .
$$
\noindent Obviously,
$$
(\tilde R +Z_R)\cdot(\tilde C +Z_C) =\tilde R \cdot \tilde C  \; ,
\eqno (5.12)
$$
\noindent and
$$
(\tilde R +Z_R)^2 =(\tilde R^2 ) +2(\tilde R \cdot Z_R) +(Z_R^2)=
(\tilde R^2) +2{\delta}_R .
\eqno (5.13)
$$
\noindent Analogously,
$(\tilde R_C +Z_C)^2 = (\tilde C^2) +2{\delta}_C $.

It remains to determine how many singular points $x\in Sing\, X $ lie
on $R $, respectively, on $C $. From 5.2 it follows that over each ss-point
on $R$ there lie 2, and on $C$ there lie $(N_2-2)$ singular points. There are no
other singular points on $R$. There are singular points on $C$ of the
following type:
over a sp-point of type $A_1$ there are $(N_2-4)$ singular points of type
$A_1 $, over a sp-point of type $A_2$ there are $(N_2-3)$ singular points of
type $A_2 $. We obtain
$$
{\delta}_R =2({\delta}_0 +n_{ss} +c_{ss}) =2({\delta}_1 -n_{sp} -c_{sp}) ,
\eqno (5.14)
$$
$$
{\delta}_C=(N_2-2)({\delta}_0 +n_{ss} +c_{ss}) +(N_2 -4)n_{sp} +(N_2 -3)c_{sp}= $$
$$
= (N_2-2){\delta}_1 -2n_{sp} - c_{sp} .
$$
Substituting $({\tilde R}^2)$ from (5.10) and ${\delta}_R$ from (5.14)
to (5.13),  we obtain
$$
(\tilde R +Z_R)^2 =2({\tilde R_1}^2) + 2n_{sp} +2c_{sp} - c_{pp}+
4({\delta}_1 -n_{sp} -c_{sp}) =
$$
$$
= 2(({\tilde R_1}^2) + 2{\delta}_1) - 2n_{sp} -2c_{sp} - c_{pp} .
$$
\noindent Analolously we find
$$
(\tilde R_C +Z_C)^2 =
(N_2-2)({\tilde R_1}^2) +
2n_{sp} - c_{pp} +
2(N_2-2){\delta}_1 -4n_{sp} - 2c_{sp} =
$$
$$
= (N_2-2)(({\tilde R_1}^2) +2{\delta}_1 ) -
2n_{sp} -2c_{sp}  - c_{pp} .
$$

Set
$$
2n_{sp} + 2c_{sp} + c_{pp} = {\iota}_1 ,
\eqno (5.15)
$$
\noindent and let
$$
g_1 =p_a(R_1) = g + {\delta}_1
\eqno (5.16)
$$
\noindent be the arithmetic genus of the curve $R_1 $. Since by (5.6)
$({\tilde R_1}^2) +2{\delta}_1 = 3\bar d +g-1 +{\delta}_1 =3\bar d +g_1 -1 $,
finally we obtain:
$$
(\tilde R +Z_R)^2 = 2(3\bar d +g_1 -1) -{\iota}_1  , \;
(\tilde C +Z_C)^2 = (N_2 -2)(3\bar d +g_1 -1) - {\iota}_1  , \;
(\tilde R +Z_R)\cdot (\tilde C +Z_C) = {\iota}_1  .
\eqno (5.17)
$$

\vspace{0.1cm}

{\bf 5.7.} {\it The self-intersection number of the divisor $\tilde R +Z_R$
is positive.}

\begin{lem}
$$
(\tilde R +Z_R)^2 >0 .
\eqno (5.18)
$$
\end{lem}

\dvo
Recall that $2\bar d = d = deg\;B $, and
${\delta}_1 ={\delta}_0 +n_{sp} +c_{sp} $. Therefore,
$$
\begin{array}{rcl}
(\tilde R +Z_R)^2 & = &2(3\bar d +g-1 +{\delta}_1) -2n_{sp} -2c_{sp}  - c_{pp} =\\
                  & = & d+(2d+2g-2) +2{\delta}_0 -c_{pp} .
\end{array}
\eqno (5.19)
$$

Now we apply the Hurwitz formula for a generic projection
$\varphi : B \to {\PP}^1 $ of the curve $B$ from a point $P\in {\PP}^2$ onto
the line $\PP^1 $, more precisely for the covering
$\bar{\varphi} : \bar B \to {\PP}^1 $, where $\bar{\varphi} =\varphi \circ n $,
and $n: \bar B \to B $ is a normalization of the curve $B $. Obviously,
the covering $\bar{\varphi}$ is ramified at the following points. Firstly,
$\bar{\varphi}$ has a ramification of the second order at points
$\bar b \in \bar B $, which correspond to non-singular points $b\in B $,
for which the line $\overline{Pb}$ is tangent to $B $.
The number of such points is $\hat d = \deg \hat B $, where $\hat B$ is a
curve dual to $B$. Secondly,  $\bar{\varphi}$ has a ramification of
order $m_k$ at points $\bar b $, which correspond to the branches $B_k$
of the curve $B$ at the singular points $b $. Here $m_k$ is the multiplicity
(order) of the corresponding branch. Denote by
$$
\nu = \sum_k(m_k-1) ,
\eqno (5.20)
$$
where the summation runs over all branches of the curve (at singular points).
The covering $\bar{\varphi}$ is of degree $d= \deg B $. By the Hurwitz
formula we obtain
$$
2g-2 = -2d+ \hat d +\nu .
\eqno (5.21)
$$

\begin{zam}
Actually we derived one of the Pl\"{u}cker formulae
$$
\hat d = 2d +(2g-2)  - \nu
$$
for a plane curve with singularities.
\end{zam}

Obviously, the number $\nu$ for A-D-E-singularities is equal:
$$
\nu (A_{2k-1}) =\nu (D_{2k+2}) =0  , \;
\nu (A_{2k}) =\nu (D_{2k+3}) =\nu (E_{7}) = 1  ,\;
\nu (E_{6}) =\nu (E_{8}) = 2  .
$$
Therefore, for the curve $B$ the number $\nu =\nu (B)$ is equal:
$$
\nu =c+\nu' , \hphantom{a} \mbox{where} \hphantom{a}
\nu' =\sum_{k>1}a_{2k} +\sum d_{2k+3} +2e_6 +e_7 +2e_8 \, .
\eqno (5.22)
$$
Returning to the proof of the inequality, we obtain from (5.19), (5.21)
and (5.22)
$$
(\tilde R +Z_R)^2 =
d+ (\hat d +\nu ) +2{\delta}_0 -c_{pp} =
d+ \hat d  +2{\delta}_0 + \nu' + (c-c_{pp}) > 0 .
\hphantom{aaa} \rule{6pt}{6pt}
\eqno (5.23)
$$

\vspace{0.1cm}

{\bf 5.8.} {\it Conclusion of the main inequality.}
Applying the Hodge index theorem to divisors
$\tilde R +Z_R $ and $\tilde C +Z_C$ on the surface $S $, we obtain
$$
\left|
\begin{array}{cc}
2(3\bar d+g_1-1)-{\iota}_1 & {\iota}_1 \\
{\iota}_1 &   (N_2-2)(3\bar d+g_1-1)-{\iota}_1
\end{array}
\right|
\le 0 .
$$
Therefore,
$$2(N_2-2)(3\bar d+g_1-1)^2 - N_2(3\bar d+g_1-1){\iota}_1 \le 0 \, $$
or
$$
N_2[2(3\bar d+g_1-1)-{\iota}_1] \le 4(3\bar d+g_1-1) \, .
\eqno (5.24)
$$
Thus, if there are two nonequivalent generic coverings
$f_1$ and $f_2$, then
$$
N_2 \le \frac{\displaystyle 4(3\bar d+g_1-1)}{\displaystyle 2(3\bar d+g_1-1)
-{\iota}_1}.
\eqno (5.25)
$$

\section{Proof of the Chisini conjecture for pluricanonical embeddings
of surfaces of general type.}

{\bf 6.1.} {\it The numerical invariants in the case of a m-canonical
embedding.}
Let $S$ be a minimal model of a surface of general type with numerical
invariants $\left( K_S^2\right) = k$ and $e(S)=e$. Let $X$ be a canonical
model of the surface $S$, and $\pi : S \to X$ be the blowing down of
(-2)-curves. Let $f: X \to {\PP}^2$ be a generic m-canonical covering, i.e.
a generic projection onto ${\PP}^2$ of $X={\varphi}_m(S)$,
where ${\varphi}_m$ is a m-canonical map,
${\varphi}_m: S \to {\PP}^{p_m-1}$, defined by the complete linear system
$\left| mK_S\right|$, $\: p_m=\frac12 m(m-1)k +\chi (S) $.
As is well known \cite{BPV} , by a theorem of Bombieri
${\varphi}_m(S)\simeq X$ for $m\ge 5$, and ${\varphi}_m$
gives the blowing down $\pi $.

Let $B\subset {\PP}^2$ be the discriminant curve. We concerve the notations
of $\S 4 $. Then
$$
\bar L =mK_S , \; K_S\cdot Z =0 , \; \bar R =(3m+1)K_S -Z .
\eqno (6.1)
$$
By formulae (4.3), we obtain
$$
N=m^2k , \; d =m(3m+1)k.
\eqno (6.2)
$$
By formulae (4.10), we find
$$
3\bar d +p_a(R) -1 = (3m+1)^2k \, ,
\eqno (6.3)
$$
and
$$
p_a(R) -1 =\frac12 (3m+1)(3m+2)k \, .
\eqno (6.4)
$$

\vspace{0.1cm}

{\bf 6.2.} {\it Invariants of a surface and of the discriminant curve
define invariants of the covering.} Now let $S_1$ and $S_2$ be two
surfaces of general type with numerical invariants $k$ and $e$. Let
$f_i: X_i \to {\PP}^2 ,$ $i=1, 2 $, be their $m_i$-canonical coverings
having the same discriminant curve $B\subset {\bf P}^2$. Then by the
second formula of (6.2) it follows that $m_1=m_2=m $. Then also
$\deg f_1 =\deg f_2 =N $. We show that the other numerical invariants of
$f_1$ and $f_2$ are the same.

By  formula (6.4) it follows that $p_a(R_1) =p_a(R_2)$, and since
$p_a(R) =g+{\delta}_X $, we have ${\delta}_{X_1} ={\delta}_{X_2} $.

By formulae (4.14) and (4.11) it follows that the number of p-cusps $c_p$ and
the number of p-nodes $n_p$ for both coverings are the same. Then
$n_p =n_{pp} +n_{ps} =n_{pp} +n_{sp} , \;
c_p =c_{pp} +c_{ps} =c_{pp} +c_{sp} $ and, consequently, $n_{ps} =n_{sp}$
and $c_{ps} =c_{sp}$.

\vspace{0.1cm}

{\bf 6.3.} {\it The main inequality in the case of surfaces of general type.}
To prove that m-canonical projections $f_1$ and $f_2$ are equivalent, by
(5.24) it is sufficient to show that an inequality
$$
N\left( 2(3\bar d +p_a(R) -1) - {\iota} \right) >
4(3\bar d +p_a(R) -1)
$$
holds (here $R$ stands for $R_1$), or
$$
(N-2)(3d +2p_a(R) -2) - N\cdot {\iota} >0,
\eqno (6.6)
$$
where
$$
{\iota} =2n_{sp} +2c_{sp} +c_{pp} =2n_{sp} +c_{sp} +c_{p}.
\eqno (6.7)
$$

Let us obtain an estimate for the number ${\iota}.$ We can express
$c_p$ by formulae (4.14)
$$
c_p =3N +2p_a(R) -2 -e.
\eqno (6.8)
$$
To estimate $2n_{sp} +c_{sp}$ we use the Hirzebruch-Miyaoka inequality
(\cite{BPV} ,p.215): if the minimal surface of general type $S$ contains
$s$ disjoint (-2)-curves, then
$$
s \le \frac29 \left( 3e(S) -(K_S^2)\right).
\eqno (6.9)
$$
Since we can take one (-2)-curve for each of the singular points of types
$A_1$ and $A_2$ on $X$, we have
$$
n_s +c_s \le \frac29 \left( 3e -k \right) .
\eqno (6.10)
$$

\begin{zam}
Instead of the Hirzebruch-Miyaoka inequality we can use the estimate
$2n_{sp} +c_{sp} \le 2(h^{1,1} -1) = 2(e-2+4q-2p_g -1)$ and the
inequalities $p_g \ge q $, $p_g \le \frac12 (K_S^2) +2$ (the Noether's
inequality).
\end{zam}

By (6.7), (6.8) and (6.10), we obtain an estimate
$$\begin{array}{rl}
{\iota} & \le
\frac49 (3e-k) +3N +2p_a(R) -2 -e = \\
& =
\frac13 e -\frac49 k +3N +2p_a(R) -2.
\end{array}
$$
Applying the Noether's inequality (\cite{BPV} ,p.211),
$$
e \le 5k +36,
\eqno (6.11)
$$
we obtain
$$
{\iota} \le \frac{11}{9}k +12 + 3N +2p_a(R) -2.
\eqno (6.12)
$$
Combining (6.12) and (6.6), we obtain a corollary.

\begin{lem}
If the inequality
$$
3N(d-N) -6d -4(p_a(R) -1) -\left( \frac{11}{9} k +12\right) N > 0
\eqno (6.13)
$$
holds, then a generic m-canonical projection of a surface of general type $S$
with given $k$ and $e$ is unique. $\rule{6pt}{6pt}$
\end{lem}

\vspace{0.1cm}

{\bf 6.4.} {\it Proof of Theorem \ref{T3} .} Express the inequality (6.13)
in terms of $m$. Substitute $N$ and $d$ from (6.2) and $p_a(R) -1 $
from (6.4) to (6.13). We obtain
$$
3m^3(2m+1)k^2 -6m(3m+1)k -2(3m+1)(3m+2)k - (\frac{11}{9} k+12)km^2 > 0,
$$
i.e.
$$
3m^3(2m+1)k -4(3m+1)^2  - (\frac{11}{9} k+12)m^2 > 0.
$$
Dividing by $m^2 $, we obtain
$$
3m(2m+1)k - \left( \frac{11}{9} k+12\right) - 4\left( 3+\frac1m \right)^2 > 0 ,
$$
or, dividing by $k$,
$$
3m(2m+1) > \frac{11}{9}  +\frac1k \left( 12 + 4(3 + \frac1m )^2 \right) .
\eqno (6.14)
$$
The right side of inequality decreases, when $k$ and $m $ increase. This
inequality holds for all $k \in {\bf N} $, if it holds for $k=1 $.
For $k=1$ and $m=3$ the right side equals
$\frac{11}{9} +12 +4\cdot \left( \frac{10}{3} \right)^2 =
\frac{173}{3} < 9\cdot 7 =63 $.
Thus, the inequality (6.14), and, consequently, the inequality (6.6),
holds for $m \ge 3$ and for all $k $. This completes the proof of
Theorem \ref{T3} .

We can mention in addition that for $m=2$ the inequality (6.14) holds,
if $k >2$, and for $m=1$ it holds, if $k>9$.

 \bigskip
 \it
 \vskip 1.5truepc
 \begin{tabular}{ll}
 {\rm V.S.~Kulikov} & {\rm Vic.~S.~Kulikov} \\
 Chair of Mathematics, & Department of algebra, \\
 Moscow State University of Printing, & Steklov Mathematical Institute, \\
 {\rm E-mail:} {\tt valentin@masha.ips.ras.ru} $\hphantom{aaa}$
 & {\rm E-mail:} {\tt victor@olya.ips.ras.ru} \\
 \end{tabular}

\end{document}